\documentclass[fleqn]{article} 
\font\BBb=msbm10 at 12pt
\newcommand{\Bbb}[1]{\mbox{\BBb #1}}

\usepackage{amscd}
                    \usepackage{amssymb}
\newcommand{\be}{\begin{equation}}
      \newcommand{\ee}{\end{equation}}
      \newcommand{\ba}{\begin{eqnarray}}
       \newcommand{\ea}{\end{eqnarray}}
\newcommand{\ban}{\begin{eqnarray*}}
\newcommand{\ean}{\end{eqnarray*}}

\newcommand{\ra}{\rightarrow}

 \newcommand{\qed}{\hspace*{\fill}\rule{3mm}{3mm}\quad \vspace{.2cm}}
 \newcommand{\Pf}{\noindent {\bf Proof:} }
 \newcommand{\Rk}{\noindent {\bf Remark} }

\newcommand{\sect}[1]{\section{#1} \setcounter{equation}{0}}

\newtheorem{theo}{Theorem}[section]
\newtheorem{defn}{Definition}[section]

\begin{document}
\newtheorem{lem}[theo]{Lemma}
\newtheorem{prop}[theo]{Proposition}
\newtheorem{coro}[theo]{Corollary}
\newtheorem{ex}[theo]{Example}
\newtheorem{note}[theo]{Note}

\title{Universal Covers for Hausdorff Limits of Noncompact Spaces}
\author{Christina Sormani\thanks
{Partially supported by NSF Grant \# DMS-0102279}
  \and Guofang Wei\thanks {Partially supported by NSF Grant \# DMS-9971833.}}
\date{}
\maketitle

\begin{abstract}
We prove that if $Y$ is the Gromov-Hausdorff limit of
a sequence of complete manifolds, $M^n_i$, with
a uniform lower bound on Ricci curvature then $Y$ has a
universal cover.\end{abstract}

\newcommand{\inj}{\mbox{inj}}
\newcommand{\vol}{\mbox{vol}}
\newcommand{\diam}{\mbox{diam}}
\newcommand{\Ric}{\mbox{Ric}}
\newcommand{\Iso}{\mbox{Iso}}
\newcommand{\Hess}{\mbox{Hess}} 
\newcommand{\divg}{\mbox{div}}
\newcommand{\RR}{\bf{R}}

\sect{Introduction}

One of the main trends in Riemannian Geometry today is the study of
Gromov Hausdorff limits. The
starting point is Gromov's precompactness theorem. Namely
a sequence of complete Riemannian $n$-manifolds with a
uniform lower bound on their Ricci curvatures have a
converging subsequence. Moreover the limit space is a
complete length space.  That is, it is a metric space such
that between every two points there is a length minimizing
curve whose length is the distance between the two points.
See \cite{Gr}, also \cite{BBI}. To prove this theorem, the
only property of Ricci curvature that Gromov uses is the
Bishop-Gromov volume comparison theorem \cite{BiCr}\cite{Gr}
which provides an estimate on the number of disjoint small
balls which fit in a large ball.

When the
sectional curvature of the sequence is uniformly bounded from
below, the limit space is well understood. Namely, it is an
Alexandrov space with curvature bounded below and, by
the work of Perelman \cite{Pl}, it is a stratified
topological manifold that is locally contractible.

In the case when only Ricci curvature is bounded from below 
Menguy has shown the limit space
can have infinite topological type on arbitrarily small
balls even with an additional assumptions of an uniform
positive lower bound on volume and nonnegative Ricci
curvature on the $M^n_i$ \cite{Me}. This is based on an
example of Perelman \cite{Per}. In the positive direction
Cheeger and Colding [ChCo1,2,3,4] have proven a number of
breakthrough results regarding the regularity and geometric
properties of the limit spaces including the construction of
a measure which satisfies the volume comparison theorem. Of
course one can not expect regularity in
a $C^{1,\alpha}$ sense but rather a statement regarding the
tangent cones at regular points. This should be contrasted
with Anderson's theorem \cite{An} that $Y$ is a
$C^{1,\alpha}$ manifold when an uniform injectivity radius
and two sided Ricci bounds are imposed on the sequence.
Despite these positive results these limit spaces are not yet
completely understood and many questions remain.

In this paper we ask whether the limit space has a
universal cover (see Definition~\ref{univcov}]).
Note that without the Ricci bound, 
a limit space may not have a universal cover.
The Hawaii Ring (Example~\ref{hawaii}) is a limit of compact  manifolds
with no curvature bound but increasingly large fundamental groups
(see \cite[Example 2.7]{SoWei}).
In contrast
authors proved that the limit spaces of sequences of simply
connected uniformly bounded compact length spaces have
universal covers which are the spaces themselves \cite[Thm
1.5]{SoWei}. With a Ricci curvature
lower bound we answer the question affirmatively
without any assumptions on simple connectivity. Namely
we prove the following theorem:

\begin{theo} \label{mainthm}
If $(Y,p)$ is the pointed
Gromov-Hausdorff limit of a sequence of $n$-dimensional
complete Riemannian manifolds $(M_i^n,p_i)$ with
$\Ric \geq (n-1)H$,  then the universal cover of $Y$ exists.
\end{theo}

In their previous article,
the authors proved this existence theorem with an additional
assumption that
the manifolds in the sequence were compact with an uniform
upper bound on the diameter \cite[Theorem 1.1]{SoWei}. To
prove that result, we defined $\delta$-covers, covers which
unravel holes of a size greater than $\delta$ [see
Definition~\ref{defdel}]. We showed that when $\delta$ is the
injectivity radius of a manifold, the $\delta$-cover is the
universal cover of that manifold. These $\delta$-covers were
complete manifolds with the same lower bound on Ricci
curvature as the original manifold, so we were able to apply
Gromov's precompactness theorem and Cheeger-Colding's
renormalized limit measures to study their limits.  We then
proved that for a compact length space that has a universal
cover, the universal cover is always some $\delta$-cover
\cite[Proposition 3.2]{SoWei}. Ultimately we showed that
there was a $\delta$ sufficiently small such that the limit
of the $\delta$-covers was the universal cover of the limit
space of the original sequence.

We cannot hope for such a strong statement in the noncompact
case. First, even for a complete manifold the universal 
cover may not be any $\delta$-cover. This can be seen from
 Nabonnand's example, which is a complete manifold
with positive Ricci curvature that is not simply connected
and yet it is its own $\delta$-cover for all values of
$\delta$ \cite{Nab}. It has a loop which is homotopic to a
sequence of increasingly small loops that diverge to the
infinity of the manifold. See also \cite{Wei}. Secondly, the
universal cover of the limit may not come from any cover of
the sequence as the following example shows. That is, pointed
limits of complete noncompact simply connected manifolds with
nonnegative sectional curvature can converge to a cylinder.

\begin{ex} \label{cylinder} If we take $M^2$,
the half cylinder capped off by a hemisphere (and suitably
smoothed, which is simply connected and has sectional
curvature $\geq 0$), and $p_i$ a sequence of points in $M$
going to infinity, then $(M^2,p_i)$ converges to a long
cylinder, its universal cover (${\mathbb R}^2$) can't come
from the limit of any cover of $M^2$ since its only cover is
$M^2$ itself. \end{ex}

To avoid these problems, it is natural to work
locally since pointed Gromov-Hausdorff convergence is
defined locally. So we will consider covers of balls. To
study the fundamental group we have to use the
intrinsic length space metric on the balls rather than the
restricted metric so we can measure the lengths of representatives of
the fundamental group. However, there has been some inconsistancy
in the literature as to whether manifolds which converge in
the pointed Gromov-Hausdorff sense have Gromov-Hausdorff
convergence of their balls with respect to the intrinsic as
well as the restricted metrics. This is clarified in our
appendix allowing us to use the intrinsic metric as long as
the radii of the balls vary in a prescribed way.

Even so, we will not use the $\delta$-covers of balls
of the converging sequence
 because it's not clear if the
sequences of such $\delta$-covers have any converging
subsequences. In general the fundmental group of a ball in
$M$ may have exponential growth even if $M$ has $\Ric \geq
0$ and is compact.

\begin{ex}  \label{ballexp}
 $M^2$ is a 2-torus obtained by gluing the sides of
a $2$ by $2$ square, take a ball at center with radius
$1<r<\sqrt{2}$, then the ball is homotopic to the
figure eight and its fundamental group has
exponential growth. 
\end{ex}
 
In fact in Example~\ref{ballnot}, we
give a sequence of smooth $2$-dimensional manifolds $M_k$
with nonnegative sectional curvature such that the universal
cover of the balls of radius $1$ in $M_k$ have no converging
subsequence. The reason that
Gromov's precompactness theorem fails for these covers is
that they are manifolds with boundary and the order of their
volume growth diverges to infinity.

Thus, in Section~\ref{delcovsect},
we introduce the relative delta cover (see
Definition~\ref{reldelcov}): the connected lift of a ball of small
radius inside the delta cover of a concentric ball of a
larger radius.  We show that if relative delta covers
converge then the limit is almost the relative delta cover of
balls in the limit (Theorem~\ref{limdelcov}). With an uniform
Ricci curvature lower bound, we can control the volumes in a
relative delta cover on the scale of
the larger ball.  We can then bound  the number of generators
of the deck transformations of the relative delta cover and,
in turn, use this to control the volume growth of the
relative delta cover, proving that relative delta covers
satisfy the Gromov precompactness theorem even though they
are manifolds with boundary with large volume growth.
This is done in Section~\ref{riccisect}.  See
Proposition~\ref{2.12a}.

In Subsection~\ref{renormsect} we extend Cheeger-Colding's construction
of a renormalized limit measure to these limit spaces although it must
be warned that these measures only satisfy the Bishop Gromov
Volume comparison theorem locally.

One might hope to construct the universal cover of the limit
space of the original manifold by piecing together the limits of relative delta
covers where the delta varies from piece to piece.  However it is possible that
these relative delta covers are nontrivial covers of balls in the universal
cover no matter how small we choose delta and how large we choose the ball.
See Example~\ref{cooper}.  Instead we use Theorem~\ref{univstab}
from Section 2, which relates the stability of the relative
delta covers with the existence of the universal cover.  We
use the renormalized limit measure to obtain the stability of
the relative delta covers. This completes the proof of
Theorem~\ref{mainthm}.

In Section~\ref{appl} we study the properties of the
universal cover and give some applications. 
First we prove that regions in the universal cover 
are limits of relative delta covers with changing
delta and $R \to \infty$.  See Theorem~\ref{univisalim}
and Corollary~\ref{2univisalim}.
Using this we are able to show that
the global Bishop Gromov
Volume comparison does hold for the lifted renormalized limit
measure on the univeral cover of the
limit space [Theorem~\ref{VOL}].

We also prove the splitting theorem on the universal cover of
the limit space when the $M_i$ have $Ricci\ge
-(n-1)\epsilon_i$ with $\epsilon_i \to 0$
[Theorem~\ref{SPLIT}] using Cheeger-Colding's almost
splitting theorem.

These two theorems allow us to generalize various results of Milnor, Anderson
and Sormani to limit spaces $Y$,
see Corollary~\ref{cormil}, \ref{corand}, and \ref{corsor}.

The authors would like to thank Professor D. Cooper for providing us with
Example~\ref{cooper}.  We would also like to thank Professor J. Cheeger
for emails concerning the details of the construction of the renormalized
limit measure in \cite{ChCo2}.

\sect{Universal Cover and $\delta$-Covering Spaces} \label{delcovsect}

\subsection{Covers and Universal Covers}

First we recall some basic definitions.

\begin{defn} {\em
We say $\bar{X}$ is a  {\em covering space} of $X$ if there
is a continuous map $\pi: \bar{X} \to X$ such that $\forall x
\in X$ there is an open neighborhood $U$ such that
$\pi^{-1}(U)$ is a disjoint union of open subsets of
$\bar{X}$ each of which is mapped homeomorphically onto $U$
by $\pi$ (we say $U$ is evenly covered by $\pi$).   }
\end{defn}

\begin{defn}\label{univcov} \cite[pp
62,83]{Sp} {\em We say $\tilde{X}$ is a
{\em universal cover} of $X$ if $\tilde{X}$ is a cover of $X$
such that for any other cover $\bar{X}$ of $X$, there is a
commutative triangle formed by a continuous map $f:\tilde{X}
\to \bar{X}$ and the two covering projections.} \end{defn}

\begin{ex}\label{hawaii}
The Hawaii Ring (see for example Spanier \cite{Sp}), is a length space
which consists of an infinite set of rings of radii decreasing to $0$
that are all joined at a common point.  This ring has no universal
cover because the universal cover would have to have an isometric
lift of a ball about this common point.  Thus the possible universal
cover  would contain many
closed rings itself and would therefore have a nontrivial cover.
\end{ex}

Let $\cal U$ be any open covering of $Y$. For any $p \in Y$, by
\cite[Page 81]{Sp}, there is a covering space, $\tilde{Y}_{\cal U}$,
of $Y$ with covering group  $\pi_1(Y,{\cal U}, p)$, where
$\pi_1(Y,{\cal U}, p)$ is a normal subgroup of $\pi_1(Y, p)$, generated
by homotopy classes of closed paths having a representative of the form
$\alpha^{-1} \circ \beta \circ \alpha$, where $\beta$ is a closed path
lying in some element of $\cal U$ and $\alpha$ is a path from $p$ to
$\beta(0)$.

Recall a few facts about covering spaces constructed from $\cal U$ \cite[Page
81]{Sp}.  If $\cal V$ is an open covering of $Y$ that refines $\cal U$, then
$\pi_1(Y,{\cal V}, p) \subset \pi_1(Y,{\cal U}, p)$, or $\tilde{Y}_{\cal V}$
covers $\tilde{Y}_{\cal U}$.
If $\pi: \bar{Y} \ra Y$ is a covering projection and $\cal U$
is an open covering of $Y$ such that each of its open sets is
evenly covered by $\pi$, then $\tilde{Y}_{\cal U}$ covers
$\bar Y$.

Recall also, that
given a metric space $(X,d)$, there is an {\em induced length metric},
$d_l$ \cite[Page 2]{Gr}, see also \cite{BBI}, given by
 \[
d_l (x,y) = \inf_{f} \{ L(f)| f: [0,1]\ra X \
 \mbox{is a
continuous map with} \ f(0)=x, f(1)=y \}, \]
where
\[
L(f) = \sup \sum_{i=0}^n d(f(t_i),f(t_{i+1})) \]
and the sup is taken among all finite partition of $[0,1]$,
$0 = t_0 \le t_1 \le \cdots \le t_{n+1} =1$. In general $d_l
\geq d$.

\begin{defn} {\em
A metric space $(X,d)$ is called a {\em length space} (or path
metric space as in \cite{Gr}) if $d = d_l$.}
\end{defn}

\begin{defn} \label{FD} {\em
Given a length space $(X,d)$ and a
covering map $\pi:\bar{X}\to X$, where $\bar{X}$ is
equipped with natural lifted length metric $\bar{d}$, the
{\em Dirichlet fundamental domain} of $X$ about $\bar{x}$,
for any $\bar{x} \in \bar{X}$, is defined: }
\be FD_{\bar{x}}=\{\bar{q} \in \bar{X}:
\bar{d}(\bar{x}, \bar{q})\le \bar{d}(\bar{z}, \bar{q})
\ \forall \bar{z}\in \bar{X}\,\,
s.t.\,\,\pi(\bar{z})=\pi(\bar{x})\}. \ee
 \end{defn}

If $Y$ is a length space with metric $d_Y$ and $X \subset Y$ then
we will use $d_Y$ to denote the restricted metric of $Y$ on $X$
and $d_X$ to denote the induced length metric or intrinsic metric
on $X$.

Let $X$ be a length space, denote $B(p, R)$ a closed ball
in $X$, i.e. $B(p, R)=\{ x\in X | d_X(x,p) \leq R\}$, and
$B_p(s)$ an open ball measured in various metrics which
will be stated on each occasion. Observe that:

\begin{lem} \label{distance}
The restricted
metric on $B(p,R)$ from $X$ is the same as the intrinsic metric on
$B(p,2R+\epsilon)$ restricted to $B(p,R)$ for any $\epsilon >0$.
Namely,
\be d_X(q_1, q_2)=d_{B(p,2R+\epsilon)}(q_1,q_2),
\qquad \forall q_1, q_2 \in B(p,R).
\ee
\end{lem}

Now let's recall the $\delta$-covers we introduced in \cite{SoWei}.
\begin{defn} \label{defdel}
Given $\delta>0$, the  $\delta$-cover, denoted $\tilde{Y}^\delta$,
of a length space $Y$, is defined to be $\tilde{Y}_{{\cal U}_{\delta}}$
where ${\cal U}_\delta$ is the open covering of $Y$ consisting of
all balls of radius $\delta$.

The covering group will be denoted $\pi_1(Y,\delta, p)\subset \pi_1(Y,p)$
and the group of deck transformations of $\tilde{Y}^\delta$ will
be denoted $G(Y,\delta)=\pi_1(Y, p)/\pi_1(Y,\delta, p)$.
\end{defn}

It is easy to see that a delta cover is a regular or Galois cover.
That is, the lift of any closed loop in Y is either always closed
or always open in its delta cover.

Note that $\tilde{Y}^{\delta_1}$ covers $\tilde{Y}^{\delta_2}$ when
$\delta_1 \leq \delta_2$.

In \cite[Prop. 3.2]{SoWei} we proved that if a compact length space $Y$
has a universal cover then it is a delta cover.  In fact $Y$ has a
universal cover iff the delta covers stabilize: there exists a
$\delta_0>0$ such that $\tilde{Y}^\delta=\tilde{Y}^{\delta_0}$
for all $\delta<\delta_0$ \cite[Thm 3.7]{SoWei}.

However, this is not true for a noncompact length space.
For example, a cylinder
with two cusped ends has ${\mathbb R}^2$ as a universal cover
but all of its delta covers are trivial.  Thus we work locally.
Since the $\delta$ covers of balls, $\tilde{B}(p,r)^\delta$,
and their covering groups $G(B(p,r),\delta)$ are not well
controlled in the Gromov Hausdorff sense as discussed in the
introduction, we define a relative $\delta$ cover:

\begin{defn} \label{reldelcov} {\em
Given closed balls
$B(p,r) \subset B(p,R) \subset Y$, where $r <R$, let
$\tilde{B}(p,r)^\delta$ and $\tilde{B}(p,R)^\delta$
be their $\delta$-covers.
The {\em relative $\delta$-cover}, denoted
$\tilde{B}(p,r,R)^\delta$, is a connected lift of $B(p,r)$ to
$\tilde{B}(p,R)^\delta$.  }
\end{defn}

Clearly relative $\delta$ covers are also regular.
The group of deck transformation of
$\tilde{B}(p,r,R)^\delta$ is denoted $G(p,r,R,\delta)$ and
is, in fact, the image of
$i_*: G(B(p,r), \delta) \ra G(B(p,R), \delta)$.
So we have the covering relation
$\tilde{B}(p,r)^\delta \ra \tilde{B}(p,r,R)^\delta \ra B(p,r)$.

\begin{lem}\label{localisom}
The covering map
\be
\pi^\delta (r,R): (\tilde{B}(p,r,R)^\delta,
d_{\tilde{B}(p,r,R)^\delta}) \ra (B(p,r), d_{B(p,r)})
\ee
is an isometry on balls of radius $\delta/3$.
\end{lem}

\Pf
If $\tilde{x} \in \tilde{B}(p,r,R)^\delta$, the intrinsic metric's ball
$B_{\tilde{x}}(s), d_{\tilde{B}(p,r,R)^\delta}$ is a subset
of the restricted ball $B_{\tilde{x}}(s), d_{\tilde{B}(p,R)^\delta}$.
Thus, for all $s<\delta$ it is mapped homeomorphically onto its image,
$U_x(s)$, in $B(p,r)$
under the map $\pi^\delta(r,R)$ (which agrees with $\pi^\delta(R)$ as
a homeomorphism).

So if $q_1,q_2 \in U_x(\delta/3)$, they lift to
unique $\tilde{q}_1, \tilde{q}_2 \in
(B_{\tilde{x}}(\delta/3), d_{\tilde{B}(p,r,R)^\delta})$.
In fact the lifts are unique in
$(B_{\tilde{x}}(\delta), d_{\tilde{B}(p,r,R)^\delta})$.
Furthermore, since
$\pi^\delta(r,R)$ is distance decreasing,
\be
d_{B(p,r)}(q_1,q_2)\le
d_{\tilde{B}(p,r,R)^\delta}(\tilde{q}_1,\tilde{q}_2) < 2(\delta/3).
\ee

Thus there is a curve $C\subset B(p,r)$ joining $q_1$ to $q_2$ of
length $L(C)=d_{B(p,r)}(q_1,q_2)<2\delta/3$.  The curve lifts to
$\tilde{C}\subset \tilde{B}(p,r,R)^\delta$ starting from
$\tilde{q}_1$ and remaining in
$B_{\tilde{q}_1}(2\delta/3)\subset B_{\tilde{x}}(\delta),
d_{\tilde{B}(p,r,R)^\delta}$.  Since the end point of
$\tilde{C}$ is a lift of $q_2$ and lifts are unique in
$B_{\tilde{x}}(\delta)$, $\tilde{C}$ joins $\tilde{q_1}$ to
$\tilde{q}_2$.  Thus \be
d_{\tilde{B}(p,r,R)^\delta}(\tilde{q}_1,\tilde{q}_2)
=L(\tilde{C})=L(C)=d_{B(p,r)}(q_1,q_2). \ee

This implies that $U_x(\delta/3)=(B_x(\delta/3), d_{B(p,r)})$
and $\pi^\delta (r,R)$ is an isometry on balls of radius
$\delta/3$. \qed

These relative delta covers provide us with a means of constructing
a universal cover.

\begin{lem}  \label{univ}
Let $Y$ be a length space. If for all $x\in Y$ there exists open neighborhood of
$x$, $U_x$, which lifts homeomorphically to all covers of $Y$,
then the cover $\tilde{Y}_{\cal U}$ created using the open
sets $U_x$ is the universal cover of $Y$.
\end{lem}
\Pf Let $\pi: \bar{Y} \ra Y$ be any covering projection.  Then by the definition
of $\cal U$, $\cal U$ is a covering
of $Y$ by open sets each evenly covered by $\pi$.  Thus,
$\tilde{Y}_{\cal U}$ covers $\bar Y$.
\qed

\begin{theo} \label{univstab}
For a length space $Y$, its universal cover
$\tilde{Y}$ exists if there is $y\in Y$ such that for
all $r>0$, there exists $R \ge r$, such that $\tilde{B}(y,r,R)^\delta$
stabilizes for all $\delta$ sufficiently small.
\end{theo}

\Rk  Note there are $Y$ such
that the universal cover of $Y$ exist but not for balls in $Y$. Such
an example can be found in \cite[Example 2.6]{SoWei}.  This example
is a compact length space which is its own universal cover but has balls
that are homotopic to the Hawaii Ring. That's one reason we use
relative delta covers instead of just delta covers of
balls.
\newline

\Pf We will construct $\tilde{Y}$ using Lemma~\ref{univ}. For any $x \in Y$,
 there exists $r$ such that $x\in B_y(r/10)$.
 Since $\tilde{B}(y,r,R)^\delta$ stabilizes, there exist $\delta_{r,R}>0$
depending on $r$ and $R$ and $\delta_x<r/10$ depending on $x,r$ and $R$
such that $B_x(\delta_x)$ lifts
homeomorphically to $\tilde{B}(y,r,R)^\delta$ for all $0<\delta <\delta_0$,
where $\delta_0=\min\{\delta_x, \delta_{r,R}\}$.
Note that by keeping $\delta_x<r/10$ we avoid having to choose a metric
and stay clear of the boundary.

Now $\tilde{B}(y,R)^\delta$ is just the disjoint union of
some copies of the relative delta cover,
$\tilde{B}(y,r,R)^\delta$, so $B_x(\delta_x)$ lifts
homeomorphically to $\tilde{B}(y,R)^\delta$ as well.

Let $\pi : \bar{Y} \ra Y$ be any cover of $Y$. Then
any connected component $\bar{B}$ of $\pi^{-1}(B(y,R))\subset
\bar{Y}$ is a covering of $B(y,R)$.
We need only show that $B_x(\delta_x)$ lifts
homeomorphically to $\bar{B}$ and thus $\bar{Y}$.

Since $B(y,R)$ is compact, we can apply the proof of \cite[Thm 3.7]{SoWei}
to say
there exists a $\delta_1>0$ such that
$\tilde{B}(y,R)^{\delta}$ covers $\bar{B}$ for all $\delta<\delta_1$.
If we take $\delta<\min\{\delta_1, \delta_0\}$ then
$B_x(\delta_x)$ lifts
homeomorphically to $\tilde{B}(y,R)^\delta$.  Thus it projects
down homeomorphically to $\bar{B}$ as well.
 By Lemma~\ref{univ}, the universal cover of $Y$ exists.
\qed

\subsection{Covers and Convergence}

Recall that $G(p,r,R,\delta)$ is the group of deck transforms
of the relative $\delta$-cover $\tilde{B}(p,r,R)^\delta$ defined
in Defn~\ref{reldelcov}.
Note that
$G(p,r,R, \delta)$ can be represented as equivalence classes
of loops based at $p$ in the small ball, $B(p,r)$, 
where $\gamma_1$ is equivalent
to $\gamma_2$ if $\gamma_2^{-1}\circ \gamma_1$ is homotopic
in the large ball, $B(p, R)$, to a loop composed of elements
of the form   $(\alpha * \beta) * \alpha^{-1}$, where $\beta$
is a closed path lying in a ball of radius $\delta$ and
$\alpha$ is a path from $p$ to $\beta(0)$.

\begin{defn} \label{defnlength}
For any $g \in G(p,r,R,\delta)$, we can define the (translative)
$\delta$-length of $g$,
\be
l(g,r,R, \delta)= \min_{q\in \tilde{B}(p,r,R)^\delta}
d_{\tilde{B}(p,r,R)^\delta}(q, g(q)).
\ee
\end{defn}

We have the following basic properties for $\delta$-length.

\begin{lem} \label{lemdel}
For all nontrivial $g \in G(p,r,R,\delta)$, the $\delta$-length of $g$,
\be
l(g,r,R, \delta) \geq \delta.  \label{lemdeleq}
\ee
For all $\delta_1 \leq \delta_2$ we have
\be
l(g,r,R,\delta_1) \geq l(g,r,R,\delta_2). \label{lemdeleq2}
\ee
\end{lem}

\Pf  Since $G(p,r,R,\delta)$ is a subgroup of $G(B(p,R),\delta)$, $g$ is also
nontrivial in $G(B(p,R),\delta)$. By \cite[Lemma 3.1]{SoWei},
\[
l(g,R,\delta)
= \min_{q \in \tilde{B}(p,R)} d_{\tilde{B}(p,R)} (q,g(q)) \geq
\delta.
\]

Clearly $l(g,r,R, \delta) \geq l(g,R,\delta)$. So $l(g,r,R, \delta) \geq
\delta$.

(\ref{lemdeleq2}) follows from that $\tilde{B}(p,r,R)^{\delta_1}$ covers
$\tilde{B}(p,r,R)^{\delta_2}$ and it's distance
nonincreasing. \qed

We now study the relationship between the relative $\delta$-covers
of two distinct balls which are close in the Gromov-Hausdorff sense,
extending Theorem 3.4 in \cite{SoWei}.  There is
some difficulty involving the use of restricted versus intrinsic length
metrics when taking Gromov-Hausdorff approximations.  So we ask the reader
to refer to the appendix, and, in particular, Defn~\ref{HausApprx} at this
time.

\begin{theo}  \label{theodelsurj}
Let $B(p_i,r_i) \subset B(p_i,R_i)\subset Y_i, i=1,2$ be
 balls each with intrinsic metrics. If there is a pointed $\epsilon$-Hausdorff
approximation $f:
B(p_1,R_1) \ra B(p_2,R_2)$ such that its restriction also gives a pointed
$\epsilon$-Hausdorff approximation from $B(p_1,r_1)$ to $B(p_2,r_2)$,
 then there is a surjective homomorphism,
$\Phi: G(p_1,r_1,R_1,\delta_1) \to G(p_2,r_2,R_2,\delta_2)$
for any $\delta_1> 10 \epsilon$ and $\delta_2>\delta_1+10\epsilon$.
\end{theo}

Note that Theorem~\ref{theodelsurj} combined with Corollary~\ref{R-r} of the
appendix gives us the following.

\begin{coro}  \label{GG}
Suppose $(M_i,p_i)$ converges to $(Y,y)$ in the pointed
Gromov Hausdorff topology, then for any $r<R$,
$\delta_1<\delta_2$, there exist
sequences $r_i\to r$ and $R_i\to R$ such that $B(p_i,r_i)$ and $B(p_i,R_i)$
converge to $B(p,r)$ and $B(p,R)$ with respect to intrinsic metrics,  and a
number $N$ sufficiently large depending upon $r,R,\delta_2$ and $\delta_1$ such
that $\forall i\ge N$ there is
 a surjective map $\Phi_i: G(p_i,r_i,R_i,\delta_1)  \to G(y,r,R,\delta_2)$.
 \end{coro}

\noindent {\bf Proof of Theorem~\ref{theodelsurj}}:
For a closed curve $\gamma: [0,1] \rightarrow B(p_1,R_1)$ with $\gamma(0)=
\gamma(1) =p_1$, construct a $5\epsilon$-partition of $\gamma$ as follows. On
$\Gamma := \gamma ([0,1])$
choose a partition $0 =t_0 \leq t_1 \leq \cdots
\leq t_m = 1$ such that for each $\gamma_i := \gamma
|_{[t_i,t_{i+1}]}, i= 0, \cdots, m-1$, one has $L(\gamma_i) <
5\epsilon$. Let $x_i = \gamma(t_i)$, and $\{x_0,
\cdots, x_m\}$ is called a $5\epsilon$-partition of $\gamma$.

For each $x_i$, we set $y_m = y_0=p_2$ and
$y_i =f(x_i), i=1,\cdots,m-1$. If $y_i, y_{i+1}$ are both in
$B(p_2,r_2)$, connect them by a minimal length curve,
$\bar{\gamma}_i$ in $B(p_2,r_2)$, otherwise connect them with a minimal length
curve, $\bar{\gamma}_i$ in $B(p_2,R_2)$. This yields a closed curve
$\bar{\gamma}$ in $B(p_2,R_2)$ based at $p_2$ consisting of $m$ minimizing
segments each having length $\le 6\epsilon$.
This construction guarantees that if $\gamma \in B(p_1,r_1)$
then $\bar{\gamma}\in B(p_2,r_2)$.

Any $\alpha \in G(p_1,r_1,R_1,\delta_1)$ 
can be represented by some rectifiable closed curve
$\gamma$ in $B(p_1,r_1)$, so we can hope to define
\[
\Phi (\alpha) = \Phi ([\gamma]) := [\bar{\gamma}] \in
G(p_2,r_2,R_2,\delta_2).
\]
First we need to verify that $\Phi$ doesn't depend on
the choice of $\gamma$ such that $[\gamma]=\alpha$.  

Using the facts that $18\epsilon<\delta_2$,
one easily see that $[\bar{\gamma}]$ doesn't depend on the choice of
minimizing curves $\bar{\gamma}_i$,
nor on the special partition $\{x_1, \cdots, x_m \}$ of $\gamma ([0,1])$.
Moreover using additionally the uniform continuity of a
homotopy one can see that $[\bar{\gamma}]$ only depends on the homotopy
class of $\gamma$ in $\pi_1(B(p_1,r_1), p_1)$.

It thus also easy to check that $\Phi$ is a homomorphism from
$\pi_1(B(p_1,r_1), p_1)$ to $G(p_2,r_2,R_2,\delta_2)$.
However $\alpha \in G(p_1,r_1,R_1,\delta_1)$ not $\pi_1(B(p_1,r_1), p_1)$.

Suppose $\gamma_1$ and $\gamma_2$ are both representatives of $\alpha
\in G(p_1,r_1,R_1,\delta_1)$.  Then $\gamma_1 *\gamma_2^{-1}$ is, in
$B(p_1,R_1)$,
homotopic to a loop $\gamma_3$ generated by loops of the
form  $\alpha * \beta * \alpha^{-1}$, where $\beta$ is a closed path lying
in a ball of radius $\delta_1$ and $\alpha$ is a path from
$p_1$ to $\beta(0)$.
So $[\bar{\gamma_1}]=[\bar{\gamma_3}] * [\bar{\gamma_2}]$ and
we need only show that $[\bar{\gamma_3}]$ is trivial in
$G(p_2,r_2,R_2,\delta_2)$.

In fact $\bar{\gamma_3}$ can be chosen  as follows.
The $y_i$'s corresponding to the $x_i$'s from the $\beta$ segments of
$\gamma_3$  are all within $\delta_1+\epsilon$ of a common point
and the minimal geodesics between them are within
$\delta_1+(1 +6/2) \epsilon < \delta_2$.  Furthermore, the  $y_i$'s
corresponding to the $x_i$'s from the $\alpha$ and $\alpha^{-1}$ segments
of the curve can be chosen to correspond.  Thus  $\bar{\gamma_3}$
is generated by loops of the
form  $\alpha * \beta * \alpha^{-1}$ lying in $B(p_2,R_2)$, where $\beta$ is a
closed path lying
in a ball of radius $\delta_2$ and $\alpha$ is a path from $p_2$ to $\beta(0)$.
So it is trivial.

Last, we need to show that $\Phi$ is onto. If
$\bar{\alpha} \in G(p_2,r_2,R_2,\delta_2)$, it can be represented by some
rectifiable closed curve $\sigma$ in $B(p_2,r_2)$ based at $p_2$.  Choose an
$\epsilon$-partition
$\{y_0, \cdots, y_m\}$ of $\sigma$. Since $f: B(p_1,r_1) \ra B(p_2,r_2)$ is an $\epsilon$-Hausdorff
approximation, there are $x_i \in B(p_1,r_1)$, $y_i'=f(x_i) \in B(p_2,r_2)$
where $y_0' = y_m' =p_2$, $x_0 = x_m =p_1$ and 
$d_{B(p_2,r_2)}(y_i,y_i') \le \epsilon$.
Connect $y_i', y_{i+1}'$ with a length minimizing curve in $B(p_2,r_2)$; this
yields a piecewise length minimizing closed curve $\sigma'$ in $B(p_2,r_2)$
based at $p_2$, each segment has length $\le 3\epsilon$. So $[\sigma'] =
[\sigma]$ in $G(p_2,r_2,R_2,\delta_2)$. Now connect $x_i, x_{i+1}$ by length
minimizing
curves in $B(p_1,r_1)$ this yields a piecewise length minimizing
$\gamma: [0,1] \rightarrow B(p_1, r_1)$ with base point $p_1$, each segment has
length $\le 4\epsilon$. So the curve $\gamma$ allows a
$5\epsilon$-partition and $[{\gamma}]\in G(p_1, r_1,R_1,\delta_1)$.
By the construction, $\Phi([\gamma])=\bar{\alpha}$.

Therefore $\Phi$ is surjective.
\qed

We prove a local relative version of Theorems 3.6 in \cite{SoWei}. Namely the
Gromov-Hausdorff limit of relative delta covers is close to being a relative
delta cover of the limit.

\begin{theo}  \label{limdelcov}
Suppose $(M_i,p_i)$ converges to $(Y,y)$ in the pointed
Gromov-Hausdorff topology.  Fix any
$R>3r>0$, and let $r_i$, $R_i$ be chosen such that
$B(p_i,r_i) \to B(y,r)$ and $B(p_i,R_i)\to B(y,R)$
with intrinisic metrics.
If $(\tilde{B}(p_i,r_i,R_i)^\delta,
d_{\tilde{B}(p_i,r_i,R_i)^\delta}, \tilde{p}_i)$ converges in
the pointed Gromov-Hausdorff metric to some space which we
will denote, $(B(p,r,R)^\delta, d_{B(p,r,R)^\delta},
\tilde{p})$, then $B(p,r,R)^\delta$ is a covering space of
$B(p,r)$. Furthermore, for all $\delta_2>\delta>\delta_1$,
$\tilde{B}(p,r,R)^{\delta_1}$ covers $B(p,r,R)^\delta$ and
$B(p,r,R)^\delta$ covers $\tilde{B}(p,r,R)^{\delta_2}$.
Namely we have covering
projections
\be  \label{coverproj}
\tilde{B}(p,r,R)^{\delta_1} \ra B(p,r,R)^\delta
\ra \tilde{B}(p,r,R)^{\delta_2}
\ra B(p,r).
\ee
\end{theo}

\Pf By Corollary~\ref{R-r}, there exist
sequences $r_i\to r$, $R_i\to R$,
and maps $f_i$ such that for all $\delta>0$ there
exists $N_\delta(r,R)$ such that for all $i \ge N_\delta(r,R)$ the maps
$f_i: B(y,R)\to B(p_i,R_i)$
are $\delta$-Hausdorff approximations with respect to
the intrinsic distances and their restrictions  $f_i:B(y,r)\to B(p_i,r_i)$ are
also $\delta$ Hausdorff approximation with respect to the intrinsic distances
on these smaller balls.
 So $B(p_i,r_i)$ and $B(p_i,R_i)$ converge to
$B(p,r)$ and $B(p,R)$ with respect to intrinsic metrics.

Let $\pi^\delta (r_i,R_i): (\tilde{B}(p_i,r_i,R_i)^\delta,
d_{\tilde{B}(p_i,r_i,R_i)^\delta}) \rightarrow
(B(p_i, r_i), d_{B(p_i, r_i)})$ be the covering map. It is
distance nonincreasing by construction. After
possibly passing to a subsequence it follows from a generalized version of the
Arzela-Ascoli theorem (see \cite{GP}, also \cite[Page 279,
Lemma 1.8]{Pe2}) that $\pi^\delta (r_i,R_i)$ will
converge to a distance nonincreasing map $\pi^\delta (r,R):
B(p,r,R)^\delta \rightarrow B(p,r)$.

By Lemma~\ref{localisom} the covering map
\be
\pi^\delta (r_i,R_i): (\tilde{B}(p_i,r_i,R_i), d_{\tilde{B}(p_i,r_i,R_i)}) \ra
(B(p_i,r_i), d_{B(p_i,r_i)})
\ee
 is an isometry on balls of radius $\delta/3$.
So the limit projection $\pi^\delta (r,R): B(p,r,R)^\delta
\rightarrow B(p,r)$ is also an isometry on balls of radius $\delta/3$
and $B(p,r,R)^\delta$ is a covering space of $B(p,r)$.

By the Unique Lifting Theorem \cite[Lemma 3.1, Page 123]{Ma}
if $\tilde{Y}_1$ and $\tilde{Y}_2$ are covers of $Y$, then $\tilde{Y}_1$
covers $\tilde{Y}_2$ if every closed curve in $Y$ which lifts to a closed
curve in $\tilde{Y}_1$ also lifts to a closed curve in $\tilde{Y}_2$.

Since $\pi^\delta (r,R)$ is an isometry on balls of radius $\delta/3$ we
 have the covering projections
\be \label{uppercover}
\pi:\tilde{B}(p,r)^{\delta/3} \ra B(p,r,R)^\delta \ra B(p,r).
\ee
To show the first projection in (\ref{coverproj}), it's enough to show that for
$\delta/3 < \delta_1 <\delta$, we have $\tilde{B}(p,r,R)^{\delta_1}$ covers
$B(p,r,R)^\delta$.  We will use (\ref{uppercover})
to view $\tilde{B}(y,r,R)^{\delta_1}$ from above.

We can look at $\tilde{B}(p,r,R)^{\delta_1}$ as
$\tilde{B}(p,r)^{\delta/3}/\sim$
where $a \sim b$ iff there is a curve, $\tilde{C}$ from $a$ to $b$
which projects to $C$ in $B(p,r)$ that is homotopic to a
combination of $\alpha\circ\beta\circ\alpha^{-1}$ in $B(p,R)$
where $\beta\subset B_q(\delta_1)$.

From Theorem~\ref{theodelsurj}, there are homomorphisms $\phi_i$ from
$G(p,r,R,\delta_1)$ to $G(p_i,r_i,R_i,\delta)$ for all $i$ large. Let $C_i =
\phi_i ([C])\in B(p_i,r_i)$. Then $C_i$ converges to $C$ and they lift as
closed curves $\tilde{C}_i$ to $\tilde{B}(p_i,r_i,R_i)^{\delta}$. Then we can
look at the limit $\tilde{C}_\infty \in B(p,r,R)^{\delta}$ which is also a
closed
curve and is in fact the lift of $C$ in $B(p,r,R)^{\delta}$. Thus
$\pi(a)=\pi(b)$ where
$\pi$ is the covering map defined in (\ref{uppercover}).  This allows
us to define a map $\pi_*$ from equivalence classes of points
in $\tilde{B}(p,r)^{\delta/3}$ to $B(p,r,R)^{\delta}$.  That is
$\pi_*:\tilde{B}(p,r,R)^{\delta_1} \ra B(p,r,R)^{\delta}$.

Since $\pi = \pi_2 \cdot \pi_*$, where $\pi_2$ is the natual covering projection
from $\tilde{B}(p,r)^{\delta/3}$ to $\tilde{B}(p,r,R)^{\delta_1}$, and $\pi$ is
a covering projection, we can show
 $\pi_*$ commutes with natural covering projections
$B(p,r,R)^{\delta}$ to $B(p,r)$ and
$\tilde{B}(p,r,R)^{\delta_1}$ to $B(p,r)$. Therefore $\pi_*$
is a covering map by \cite[Page 131, Lemma 6.7]{Ma}.

We now prove the other part of covering maps in (\ref{coverproj}):
\be
B(p,r,R)^\delta \ra \tilde{B}(p,r,R)^{\delta_2} \ra B(p,r)\
\forall \delta_2>\delta.
\ee

Suppose not there is a $\delta_2$ for which it is not a covering map.
Then there is a closed curve $C$ in
$B(p,r)$ whose lift to $B(p,r,R)^\delta$ is closed but
whose lift to $\tilde{B}(p,r,R)^{\delta_2}$ is not a closed loop.

Since the lift of $C$ in $\tilde{B}(p,r,R)^{\delta_2}$ is not closed,
$\Phi_{\delta_2}([C]) \in G(p,r,R,\delta_2)$ is nontrivial.
Using Corollary~\ref{GG},
we can find $N$ sufficiently large so that
$\Phi_i:G(p_i,r_i,R_i,\delta)  \to G(p,r,R,\delta_2)$
is surjective. In particular we can find curves $C_i$ which converge
to $C$ in the Gromov-Hausdorff sense, such that $\Phi_i([C_i])=[C]$.
Since, $[C_i]$ are nontrivial their lifts to
$\tilde{B}(p_i,r_i,R_i)^\delta$
run between points $\tilde{C}_i(0)\neq\tilde{C}_i(1)$.

Furthermore, by
Lemma~\ref{lemdel}
\be
d_{\tilde{B}(p_i,r_i,R_i)^\delta}(\tilde{C}_i(0),\tilde{C}_i(1))
=d_{\tilde{B}(p_i,r_i,R_i)^\delta}(\tilde{C}_i(0),[C_i]\tilde{C}_i(0))
\ge l([C_i],r_i,R_i,\delta)\ge \delta.
\ee
In the limit, the lifted curves $\tilde{C_i}$ converge to the
lift of the limit of the curves, $\tilde{C}$ in $B(p,r,R)^\delta$  and
$$
d_{B(p,r,R)^\delta}(\tilde{C}(0), \tilde{C}(1))\ge \delta.
$$
This implies that $\tilde{C}$ is not closed and we have a contradiction.
\qed

\sect{Relative $\delta$-covers with $\Ric \ge (n-1)H$}
   \label{riccisect}

\subsection{Gromov's Precompactness Extended}

In order to apply Theorem~\ref{limdelcov}, we need to prove that
sequences of relative delta covers have Gromov-Hausdorff limits
even though they are manifolds with boundary
(Proposition~\ref{2.12a}). We show this can be done in the case when
the balls are in manifolds
with lower bounds on Ricci curvature.  Recall that Gromov's
Precompactness Theorem states that if $M^n_i$ are a sequence
of complete $n$-dimensional manifolds
with $\Ric \ge (n-1)H$, then a subsequence of
 $M^n_i$ converge to a complete length space $Y$ \cite{Gr}.
There are two equivalent ways to describe this convergence of
noncompact spaces as can be seen in Defn~\ref{GrRest},
Defn~\ref{GrInt} and Lemma~\ref{grdefsame}
of the appendix.  

We now state our precompactness result.

\begin{prop} \label{2.12a}
Let $M^n_i$ be complete $n$-dimensional Riemannian manifolds
with $\Ric \ge (n-1)H$ that converges to a length space $Y$.
Fix any $R>0$ and $r \in (0,R/3)$.  Let $r_i$ and $R_i$ be
radii such that $B(p_i,r_i)$ and $B(p_i,R_i)$ in $M_i$ with
intrinsic metrics converge to $B(y,r)$ and $B(y,R)$ in $Y$
with intrinsic metrics. Then for all $\delta\in (0,r/2)$ the
relative $\delta$ covers,
$(\tilde{B}_i(p_i,r_i,R_i)^\delta, \tilde{p}_i)$, with
intrinsic metric $d_{\tilde{B}_i(p_i,r_i,R_i)^\delta}$, have
a convergent subsequence in the pointed Gromov Hausdorff
sense. \end{prop}

Note that the universal cover of the balls
themselves need not have a converging subsequence in the
pointed Gromov-Hausdorff sense
even when they are given a nonnegative sectional curvature condition
[Example~\ref{ballnot}].

To prove Proposition~\ref{2.12a} we need several lemmas. By
Gromov's Precompactness Theorem \cite[Page 280, Lemma 1.9]{Pe2} [Gr],
it is enough to bound the number of $\epsilon$-net points, the
centers of the minimal set of $\epsilon$-balls whose union covers
a ball of a given radius centered at $\tilde{p_i}$ in
$\tilde{B}_i(p_i,r_i,R_i)^\delta$, uniformly for all $i$. 
This can be done by bounding the maximum number of disjoint balls
of $\epsilon/2$ in any ball of fixed radius $B_{\tilde{p}_i}(s)$.

Using the Bishop Gromov Volume Comparison Theorem, the $\epsilon$-net
points can be easily bounded in arbitrarily large balls within complete
spaces with a lower bound on Ricci curvature.
However, $\tilde{B}(p,r,R)^\delta$ is not a complete manifold and nor is
$\tilde{B}(p,R)^\delta$.  Note that $\tilde{B}(p,R)^\delta$
could have {\em exponential}, not polynomial,
volume growth (see Example~\ref{ballexp}).  Nevertheless, we 
can prove Bishop-Gromov Volume Comparison does hold for small
balls in $\tilde{B}(p,R)^\delta$ and then use other
techniques to count larger nets.

We begin with  some illustrative examples.

\begin{ex} \label{ballnot}

Let  $Y$ be a singular
${\mathbb R}P^2$, viewed as gluing the opposite points of the
boundary of a flat disk with center $y$ and radius $1$.

 Now let $X_k$ be similar copies of singular ${\mathbb R}P^2$
by gluing the opposite points of a regular $2k$-polygon
inscribed in the unit disk.  Let $x_k$ be the center points and
$r_k = \frac{\cos \frac{\pi}{2k}+3}{4}$ so that $B(x_k,r_k)$
reaches the boundary of the $2k$ polygon but is within the unit disk.
Then the fundamental
groups of $B(x_k,r_k)$ in $X_k$ has a free group with $k$
generators.

Clearly $(X_k, x_k)$ converge to $(Y, y)$ in Gromov Hausdorff
sense but the universal covers of
$B(x_k,r_k)$, $\tilde{B}(x_k,r_k)$, do not converge in the Gromov Hausdorff
sense and nor do  they have any converging subsequence.
This can be seen because the number of disjoint balls of
radius $cos(\pi/(2k))$ in $B_{\tilde{x}_k}(2)\subset \tilde{B}(x_k,r_k)$
is at least $2k+1$ and diverges to infinity.

Now $X_k$ can be smoothed to manifolds $M_k$ which are
diffeomorphic to ${\mathbb R}P^2$ and has nonnegative sectional
curvature. To construct $M_k$, first smooth the $X_k$ at the
corner of the polygons, then curve down the boundary so that the
boundary is totally geodesic and a reflection of this piece
across the boundary is a smooth $S^2$.
Let $M_k$ be this sphere's ${\mathbb Z}_2$ quotient. If we do the
smoothing so that region of smoothing is very small then
$M_k$ also converges to $Y$ and the fundamental
group of $B(p_k,r_k')$ in $M_k$ with $r_k' = \frac{3\cos
\frac{\pi}{2k}+1}{4}$   also has a free group with k
generators. By rescaling the manifolds $M_k$ we have a sequence
of smooth 2 dimensional manifolds $M_k$ with nonnegative
sectional curvature such that the universal cover of the
balls of radius $1$ in $M_k$ have no converging subsequence.
\end{ex}

The following example shows the difference between a
delta cover and a relative delta cover, illustrating why
balls are controlled better in the latter case.

\begin{ex} \label{4torus}
Let $M=T^4$ such that the first two circles have diameter $1$
and the second two have diameter $5$ as length spaces.  Let
$R=8$ and $r=2$.  Note that $\tilde{B}(p,R)$ and
$\tilde{B}(p,r)$ both have exponential volume growth since
 $\tilde{B}(p,R)={\mathbb R}^2 \times X$ where X branches in two
dimensions causing exponential growth on top of the polynomial
growth of order 2, and
$\tilde{B}(p,r)$ also has the
branching effect in the first two dimensions but does not
have a Euclidean factor.
However the relative $\delta$ cover is under control because
the nontrivial directions of the connected lift of $B(p,r)$ to
$\tilde{B}(p,R)$ are lifted into the Euclidean factor and
thus cannot grow more than polynomially.
\end{ex}

In the
following lemmas we will omit
the subscript $i$ and just proceed to find a uniform bound
on the number of minimal set of $\epsilon$-balls whose union covers
a ball of fixed radius centered at $\tilde{p}$ in
$\tilde{B}(p,r,R)^\delta$
depending only on $n$, $r$, $R$ and $H$.

We first have a lemma concerning bounded balls
measured using the restricted metric.

\begin{lem} \label{volume}
Let $M^n$ be a manifold with $\Ric \ge (n-1)H$.
If $B(\tilde{p},s) \subset
\tilde{B}(p,R)^\delta$ is centered on a lift of $p$
and is measured with the
metric restricted from $d_{\tilde{B}(p,R)^\delta}$, then
we have
\be \label{vol1}
\frac {\vol(B(\tilde{p},s_1))}{\vol(B(\tilde{p},s_2))}
\ge V(n,H,s_1)/V(n,H,s_2) \qquad \forall s_1<s_2\le R.
\ee
Also for all $\tilde{x} \in \tilde{B}(p,R)^\delta$, we have
\be \label{vol2}
\frac {\vol(B(\tilde{x},s_1))}{\vol(B(\tilde{x},s_2))}
\ge V(n,H,s_1)/V(n,H,s_2) \qquad \forall s_1<s_2\le R-h
\ee
where $h$ is the distance from $\tilde{x}$ to the nearest
lift of $p$, and $V(n,H,r)$ is the volume of balls of radius $r$ in the model space
$M_H^n$, the $n$-dimensional Riemannian manifold with constant
sectional curvature $H$.
\end{lem}

\Pf
 Note that the proof of Bishop-Gromov Volume Comparison \cite{BiCr}\cite{Gr}
can easily be applied here as long as we have smooth
minimal geodesics, so we need only avoid the boundary of
$\tilde{B}(p,R)^\delta$. That is we can apply it to balls
$B(\tilde{x}, s)$ such that all points $\tilde{q}$ in the
ball are joined to $\tilde{x}$ by a minimal geodesic that
avoids the boundary of $\tilde{B}(p,R)^\delta$.  This works
for $\tilde{x}=\tilde{p}$ as long as $s$ is less than $R$
since there will be a curve of length $s$ joining $\tilde{p}$
to $\tilde{q}$ and this curve will project to a curve in
$B(p,R)$ starting from $p$. If the curve hits the boundary,
then its projection hits the boundary of $B(p,R)$, so it has
length $\ge R$.

For arbitrary $\tilde{x}$, we avoid the boundary by staying
inside a ball $B(g\tilde{p}, R)$ for the lift of $p$ closest
to $\tilde{x}$. \qed

This has a corollary which gives a volume estimate for bounded balls on
$\tilde{B}(p,r,R)^\delta$ with its intrinsic length metric.  We state
four versions of the estimates because they will be needed
later to construct a measure on the limit space.

\begin{coro}\label{volcor}
Let $B(\tilde{x},s)\subset \tilde{B}(p,r,R)^\delta$
be a ball of radius $s$ measured using
$d_{\tilde{B}(p,r,R)^\delta}$.

Then if $\tilde{x}=\tilde{p}$, $s_1<r$ and $s_1<s_2<R$ we
have \be \label{volcor1}
\frac {\vol(B(\tilde{p},s_1))}{\vol(B(\tilde{p},s_2))}
\ge V(n,H,s_1)/V(n,H,s_2).
\ee
Furthermore for $\tilde{x}_1, \tilde{x}_2 \in
\tilde{B}(p,r,R)^\delta$ with
$s=d_{\tilde{B}(p,r,R)^\delta}(\tilde{x}_1,\tilde{x}_2)$ we
get the following inequalities:
\be \label{chco1.2}
\frac
{\vol (B(\tilde{x}_1,r_1))}{\vol (B(\tilde{x}_2,r_2))}
\ge \frac{V(n,H,r_1)}{V(n,H,r_2+s)} \quad  \forall r_1\le r_2 \le R-r-s, \
r_1<r-d_M(\pi(\tilde{x}_1), p),
\ee

\be \label{chco1.3}
\frac{\vol (B(\tilde{x}_2,r_2))}
{\vol (B(\tilde{x}_1,r_1))}
\ge \frac{V(n,H,r_2)}{V(n,H,r_1+s)}  \quad \forall r_2\le r_1 +s\le R-r, \
r_2<r-d_M(\pi(\tilde{x}_2), p),
\ee

\be \label{chco1.4}
\frac
{\vol (B(\tilde{x}_2,r_2))}{\vol (B(\tilde{x}_1,r_1))}
\ge 1
\qquad \forall r_2\ge r_1+s.
\ee
\end{coro}

\Pf
Note that if $B(\tilde{p},s)$
is measured with $d_{\tilde{B}(p,r,R)^\delta}$ then
it is a subset of  $B(\tilde{p},s)$
with $d_{\tilde{B}(p,R)^\delta}$.
For $s_1<r-d_M(\pi(\tilde{x}),p)$, $B(\tilde{p},s_1)$ is the
same whether it is measured with
$d_{\tilde{B}(p,r,R)^\delta}$ or with
$d_{\tilde{B}(p,R)^\delta}$ because it avoids the boundary.
Thus (\ref{volcor1}), follows from Lemma~\ref{volume} which
only allows $s_2\le R$.

We now prove the equations needed to create a measure on the limit
spaces.
Before applying (\ref{vol2}), we note that $h$
does not depend on the metric which is used
to measure it, intrinsic or restricted.
Furthermore, here $h\le r$.

We get (\ref{chco1.2}) by setting $s_1=r_1<r$, $s_2=r_2+s$
and $\tilde{x}=\tilde{x}_1$.  So we need $r_2+s\le R-r$ and
$r_1<r-d_M(\pi(\tilde{x}),p)$.

We get (\ref{chco1.3}) by setting $s_1=r_2$, $s_2=r_1+s$
and $\tilde{x}=\tilde{x}_2$.  So we need $r_1+s\le R-r$ and
$r_2<r-d_M(\pi(\tilde{x}),p)$. \qed

At this point we need to control larger regions than just balls
of radius $<R$.  We do not have a volume comparison theorem on
this scale, but we can control these regions by piecing together
fundamental domains and controlling the generators $g\in
G(p,r,R,\delta)$ which map a fundamental domain based
at $\tilde{p}$ to an adjacent fundamental domain based at $g\tilde{p}$.

\begin{lem} \label{generator}
Let $M^n$ ba a Riemannian manifold with $\Ric \ge (n-1)H$.
Let FD be the Dirichlet fundamental domain of $B(p,r)$ based
at $\tilde{p}\in \tilde{B}(p,r,R)^\delta$ with $3r<R$ (see
Definition~\ref{FD}).

Then the number of $gFD\subset B(\tilde{p}, 3r)$
with the ball measured using the
restricted metric,
$d_{\tilde{B}(p,r,R)^\delta}$, is uniformly bounded by
a number  $N= N(n,H,r,\delta)$.  So is
the number of $gFD$ of the connected lift
which are adjacent to $FD$, (i.e. $gFD\cap FD\neq \emptyset$).
 \end{lem}

\Pf Note that if $gFD$ is adjacent to $FD$, then
$gFD \subset B(\tilde{p}, 3r)$ with the ball measured using the
restricted metric,
$d_{\tilde{B}(p,r,R)^\delta}$.  Furthermore $\delta$ balls
around $\tilde{p}$ and $g\tilde{p}$ are isometric and
disjoint. So we can apply the volume comparison above
[Lemma~\ref{volume}] and packing arguments to get that
\be N \leq \vol (B(\tilde{p}, 3r))/ \vol (B(\tilde{p}, \delta)) \leq
V(n,H,3r)/V(n,H, \delta).\ee
\qed

We also need the following net lifting lemma which does not require a Ricci
curvature bound.

\begin{lem} \label{netlift}
If $\cal A$ is an $\epsilon$-net points of $B(p,r-\epsilon)$
(we can assume that
${\cal A} \subset B(p,r-\epsilon)$) with
restricted metric from $d_M$
and $\epsilon < \delta$,
then the lift of $\cal A$ is a $2\epsilon$-net points of
$\tilde{B}(p, r, R)^\delta$ and its fundamental domains with respect to the
intrinsic metric
$d_{\tilde{B}(p,r,R)^\delta}$.  Furthermore the number of $2\epsilon$-balls
needed to cover each $gFD$ is less than or equal to $N card( \cal A )$ where
$N$ is defined in Lemma~\ref{generator}.
\end{lem}

\Pf Since $\epsilon<\delta$,
the covering map $\pi: \tilde{B}(p,r,R)^\delta \ra B(p,r)$
is a diffeomorphism from $U_{\tilde{q}}$ to $B_q(\epsilon),
d_M$, where $U_{\tilde{q}}$ is the connected component of a
lift of $B_q(\epsilon)$ centered at a lift $\tilde{q}$ of
$q$.

Note that $U_{\tilde{q}}$ is not necessarily isometric to nor contained in
a ball of radius $\epsilon$ measured
using $d_{\tilde{B}(p,r,R)^\delta}$.  However, it is easy to see that
\be \label{FDcap}
FD \cap \pi^{-1}(B(p,r-\epsilon))
\subset \bigcup_{q\in \cal A} \,\,
\bigcup_{g\, s.t. gFD \cap FD \neq \emptyset} U_{g\tilde{q}}
\ee
which is a union of $N card(\cal A)$ sets.  Here each $\tilde{q}$ is a
chosen lift of $q$ to $FD$.

For any $\tilde{y}\in FD\subset \tilde{B}(p,r,R)^\delta$, there is $y\in B(p,r)$
which is joined by a minimal geodesic to $p$.  Thus there is a point
$x\in B(p, r-\epsilon)$ joined to $y$ by a smooth
geodesic of length $\le\epsilon$.
This lifts upstairs and we get a point
$\tilde{x}\in FD \cap \pi^{-1}(B(p,r-\epsilon)) $
joined to $\tilde{y}$ by a smooth geodesic
contained in $\tilde{B}(p,r,R)^\delta$ of length $\le\epsilon$.

For any $\tilde x \in FD \cap \pi^{-1}(B(p,r-\epsilon))$, by (\ref{FDcap}),
there exists a $q \in \cal{A}$ and a $g$ such that
$\tilde{x}\in U_{g\tilde{q}}$.  So $\tilde{x}$ projects to a point $x$ in
$B_q(\epsilon)\subset B(p,r)$ measured using $d_M$.  Thus $x$ has a
smooth minimal geodesic
contained in $B(p,r)$ of length $<\epsilon$ joining it to $q$, and this
lifts to a smooth curve in $\tilde{B}(p,r,R)^\delta$ joining $\tilde{x}$ to
$g\tilde{q}$.  In conclusion, for all $\tilde{y}\in FD$ we have a curve of
length $<2 \epsilon$ contained
in $\tilde{B}(p,r,R)^\delta$ joining $\tilde{y}$ to a $g\tilde{q}$:
\be \label{FDsub}
FD\subset \bigcup_{q\in \cal A} \,\,
\bigcup_{g\, s.t. gFD \cap FD \neq \emptyset}
B_{g\tilde{q}}(2\epsilon), \ee
where the balls in the union are measured using the intrinsic metric
$d_{\tilde{B}(p,r,R)^\delta}$.
\qed

Now we are ready to prove Proposition~\ref{2.12a}.
\newline

\noindent {\bf{Proof of Proposition~\ref{2.12a}:}}
Let $B(\tilde{p},\bar{R})$ be any large ball in $(\tilde{B}(p,r,R)^\delta,
d_{\tilde{B}(p,r,R)^\delta})$. Note that any $\epsilon/2$-net of
$B(p,r-\epsilon)$ with
restricted metric from $d_M$ is uniformly bounded by volume comparison
\cite{Gr}.
Hence we
just need to bound the number of fundamental domains intersecting
$B(\tilde{p},\bar{R})$ so we can lift the $\epsilon/2$-net from
$B(p,r-\epsilon/2)$ and
give a bound on $\epsilon$-net points of $B(\tilde{p},\bar{R})$
using Lemma~\ref{netlift}.

Suppose FD is a fundamental domain of $B(p,r)$ based at $\tilde{p}$
and $g$FD intersects $B(\tilde{p},\bar{R})$.  Then there is a
curve of length $\le \bar{R}$ running from a point in $g$FD
back to $\tilde{p}$ which stays
in $\tilde{B}(p,r,R)^\delta$. This curve is contained in a union of $k$
fundamental domains each of which is adjacent to the next. These fundamental
domains can be called $h_jFD$ with $h_0=e$ and $h_k= g$. So we can write $g=
g_1g_2....g_k$, where $g_1 = h_1, g_2 = h_1^{-1}h_2, \dots, g_k =
h_{k-1}^{-1}h_k$. Note that each $g_j$ is a generator (that is $g_j$FD is
adjacent to FD). So if we can find a uniform bound on $k$ then the number of
fundamental domains $g$FD
intersecting the ball of radius $\bar{R}$ is bounded by the number of
words of length $k$ in $N$ generators, which is $N^k$.

By Lemma~\ref{generator}, $N \leq N(n,H,r, \delta)$.   We
also get $\bar{N}$, the number of deck transforms $g$ such that
$g \tilde{p}$ is in $\tilde{B}(p,r,R)^\delta$ and $d(g\tilde{p}, \tilde{p}) <
3r$, is $\leq N(n,H,r, \delta)$.

We will show that $k < \bar{N}([\bar{R}/r]+1)$, where $[\bar{R}/r]$ is the
integer part of $\bar{R}/r$.  If not, then look at a sequence of points $q_0
=\tilde{p}, ..., q_j \in h_jFD, ...$ all running along the curve of
length $\le
\bar{R}$.  The last point is $q_k$.  If $k \ge \bar{N}([\bar{R}/r]+1)$ then
we can look at the points $q_0, q_{\bar{N}}, q_{2\bar{N}}, q_{([\bar{R}/r]+1)
\bar{N}}$. Each pair of these points has $\bar{N}-1$ points lying between them.
Furthermore $\min d(q_{j\bar{N}}, q_{(j+1)\bar{N}}) \le
\bar{R}/([\bar{R}/r]+1)<r$,
otherwise $d(q_0, q_k) >\bar{R}$.  Thus there exist $\bar{N} +1$ points
which are within a distance $r$ from the minimizing $q_{j\bar{N}}$.
(That is, all the points up to and including $q_{(j+1)\bar{N}}$).

Now each of these $q_k \in h_kFD$, so $d(h_k\tilde{p}, q_k) <r$.  Thus
we know there exist $\bar{N} +1$ points of the form $g\tilde{p}$
which are within a distance $3r$ from $h_{j\bar{N}} \tilde{p}$ in the length
metric and thus also in the restricted metric.  Multiplying
all points by $h_{j\bar{N}}^{-1}$ we contradict the definition of $\bar{N}$.
Thus the claim which bounds $k$ is correct, and we are done.
\qed

\subsection{Renormalized Measures} \label{renormsect}

By the previous two sections, we know that if $M^n_i$ are
complete manifolds with $\Ric \ge (n-1)H$ then for any
$p_i\in M_i$ and any $R>3r>0$, there is a subsequence of the
$i$ such that $(M_i,p_i)$ converge to $(Y,y)$ and there exist
$r_i\to r $ and $R_i\to R$ such that $(B(p_i,r_i),p_i) \to
(B(y,r),y)$ and the relative delta covers $\tilde{B}(p_i,r_i,
R_i)^\delta \to B(y,r,R)^\delta$, where $B(y,r,R)^\delta$
covers $B(y,r)$ for any fixed $\delta>0$. We now construct a
renormalized limit measure on $B(y,r,R)^\delta$ similar to
the one used by Cheeger and Colding to construct a limit
measure on $Y$. In fact we prove the following more general
theorem which allows us to vary the $\delta$ in the sequence
of relative delta covers.

\begin{prop}\label{limmeas}
Let $M^n_i$ be complete manifolds with $\Ric \ge(n-1)H$, such
that $GH\lim_{i\to\infty}(M_i,p_i)=(Y,y)$. Fix $R\in
(0,\infty]$. Then for all $r\in (0,R/4)$ and $\delta_i\in
(0,r/2)$ and $r_i\to r$ such that
\be \label{limmeasa}
(B(y,r), d_{B(y,r)}) = GH\lim_{i\to \infty} (B(p_i,r_i), d_{B(p_i,r_i)}),
\ee
 the relative $\delta_i$ covers converge
to a covering space $\hat{B}(y,r,R)$ for some
$R_i \to R\in (0,\infty]$ as follows:
\be \label{limmeasb}
(\hat{B}(y,r,R), \hat{y}, d_{\hat{B}(y,r,R)}):=GH\lim_{i\to
\infty} (\tilde{B}_i(p_i,r_i,R_i)^{\delta_i}, \tilde{p}_i,
d_{\tilde{B}_i(p_i,r_i,R_i)^{\delta_i}}),
\ee
and $\pi: \hat{B}(y,r,R)\to B(y,r)$ is a limit of covering
maps $\pi_i: \tilde{B}(p_i,r_i,R_i)^{\delta_i}\to B(p_i,r_i)$
all of which are isometries on balls of some common radius
$\delta/2>0$,
then $\hat{B}(y,r,R)$ has a renormalized limit measure,
$\mu$, which is Borel regular.
This measure is a limit of measures on a subsequence, $i_j$ of the
original sequence in the following sense:  for all $
\hat{x}\in \hat{B}(y,r,R), s<r-d_Y(\pi(\hat{x}),y)$
there exists $\tilde{x}_{i_j}\in
\tilde{B}(p_{i_j},r_{i_j},R_{i_j})^{\delta_{i_j}}$ such that
\be \label{vinfeq} \mu(B(\hat{x},s))=\bar{V}_\infty
(B(\hat{x},s))=\lim_{i_j\to\infty}
\frac{Vol(B(\tilde{x}_{i_j},s))}{Vol(B(\tilde{p}_{i_j},r/10))
} . \ee In fact $\mu$ is created from its measure on these
small balls using Caratheodory's Construction.

Furthermore, we have the Bishop-Gromov Volume Comparison,
\be \label{mu1.2}
\frac
{\mu(B(\hat{x},r_1))}{\mu(B(\hat{x},r_2))}
\ge V(n,H,r_1)/V(n,H,r_2) \ \ \forall r_1\le r_2 < (R-r)/2,
\quad
r_1<r-d_Y(\pi(\hat{x}),y). \ee
Finally $\mu$ is Radon when restricted to closed sets
contained in balls about $\hat{y}$ that avoid the boundary.
\end{prop}

Note that we must assume that the limit of the relative
delta covers exists and is a cover for this result,
as can be seen in Example~\ref{ballnot}.  Note also that when $\delta_i=\delta$
then $\hat{B}(y,r,R)=B(y,r,R)^\delta$ of
Theorem~\ref{limdelcov} and all the conditions of
Proposition~\ref{limmeas} are satisfied. However, we do not
in general assume that $\delta_i$ are bounded below by some
$\delta$ just that all $\pi_i$ are isometries on
$\delta/2$ balls.

Recall that Caratheodory's Construction consists of taking a function
$\psi:F \to {\Bbb R}$ where $F$ is a collection of  sets and
then taking an infimum
as follows:
\be \label {rcarath1}
\mu_\epsilon(A)=\inf \left\{ \sum_{B\in G} \psi(B)
                 \,\,:\,\,G\subset F \cap \{B: diam(B)\le
\epsilon\}                   \textrm{ and } A \subset
\bigcup_{B\in G} B  \right\}, \ee
and let
\be \label{rcarath2}
\mu(A)=\lim_{\epsilon\to 0} \mu_\epsilon(A).
\ee
When all the members of $F$ are Borel sets, $\mu$ is a Borel
Regular measure \cite[2.10]{Fed}.  This measure is Radon when restricted to
measurable sets by Thm 13.7 of \cite{Mun} which only requires that
the members of $F$ are open sets.

Cheeger and Colding defined a function $\psi$ on all balls in
$Y$ using the relative volume comparison theorem and taking
limits of subsequences. They then quote a standard packing
argument to create an uniform approximation of the infimum in
(\ref{carath1}). We do not have control over all balls, just
those that avoid the boundary of
$\tilde{B}(p_i,r_i,R_i)^{\delta_i}$.  So in our case $F$ does
not consist of all balls.  Thus our measure only agrees with
the one defined by Cheeger-Colding on sets whose tubular
neighborhood's avoid the boundary as follows.

\begin{coro}\label{agree}
If the renormalized limit measure, $\mu_Y$, on $Y$ is defined as
in \cite{ChCo2} with respect to balls of radius $1$ and we
take $\delta<\min\{r/10, 1\}$ in Proposition~\ref{limmeas},
then
the measure $\mu_{\hat{B}(y,r,R)}$ on $\hat{B}(y,r,R)$ of
Proposition~\ref{limmeas} agrees up to a scale with $\mu_Y$
when evaluated on closed isometrically lifted sets $S$
contained in closed balls that avoid the boundary of
$B(y,r)$. Namely
 \be
\mu_{\hat{B}(y,r,R)}(S)=\lambda \,\, \mu_Y(\pi(S)),
\ee
 where $
\lambda =\lim_{i\to\infty} \frac
{\vol(B(p_i,1))}{\vol(B(\tilde{p}_i,r/10))} \in \left[
\frac{V(n,H,\delta/2)}{V(n,H,r/10)},\frac{V(n,H,1)}
{V(n,H,\delta/2)} \right]$. \end{coro}

This corollary will follow from the proof of Proposition~\ref{limmeas}
because the definition of $\bar{V}_\infty$
is exactly as in \cite{ChCo2} and the Caratheodry
constructions used in their paper and here will agree on sets
which avoid the boundary. The upper bound for $\lambda$ is
found using Bishop Gromov on $B(p_i,1)$ and the fact that
$B(p_i,\delta/2)$ lifts isometrically and that
$B(\tilde{p}_i, \delta/2)\subset B(\tilde{p}_i,r/10)$.  The
lower bound for $\lambda$ is found using Bishop Gromov on
$B(\tilde{p}_i,r/10)$ and the fact that
$B(\tilde{p}_i,\delta/2)$ is mapped isometrically to
$B({p}_i, \delta/2)\subset B({p}_i,1)$.

The next corollary concerns both our measure and the
one defined by Cheeger-Colding.
That is $X_i$ can be taken to be complete manifolds converging to
$X$ or $X_i$ can be $\tilde{B}(y_i,r_i,R_i)^{\delta_i}$
converging to $X=\hat{B}(y,r,R)$.  In the former case
$int(X)=X$. The corollary will be proven at the end of this
subsection.

\begin{coro} \label{TKi}
Suppose $X_i$ converges to
$X$ with the renormalized measure convergence
defined above.
If $K_i\subset X_i$ converges to
a compact set $K\in int(X)$ as subsets of $X_i$ and $X$,
then there exists a sequence $\epsilon_i>0$ converging to $0$
such that
\be
\lim_{i\to\infty}
\frac{\vol(T_{\epsilon_i}(K_i))}{\vol(B(\tilde{p}_i,r/10))}=\mu (K). \ee
\end{coro}

We need the following packing lemma before we can
prove Proposition~\ref{limmeas}.

We first define special compact subsets
of covers of $B(p,r)$ which avoid the boundary by a definite amount.
Let $\pi:\tilde{B} \to B(p,r)$ be any cover and let
$R_0<r$.  Then we define \be \label{KR0}
K_{R_0}:=\{\tilde{x}\in \tilde{B}:
d_{B(p,r)}(\pi(\tilde{x}),p)\le R_0\}, \qquad
i\tilde{B}:=\bigcup_{R_0<r} K_{R_0}.
\ee
Note that
\be \label{KR0prop}
\pi(i\tilde{B})=B_p(r) \textrm{ and } Cl(i\tilde{B})=\tilde{B}.
\ee

\begin{lem} \label{packing}
Given $R > 2r$, $R_0<r$, $B(\tilde{x},s)\subset
\tilde{B}(p_i,r_i,R_i)^\delta$ (the subscript $i$ will be
omitted below in this lemma) such that $2s<R-r$. For any
compact set $K\subset B(\tilde{x},s)\cap K_{R_0}$,
$\epsilon\in (0, \min \{r-R_0, r/2\})$ there exists
$\lambda(\epsilon,s,n,H)$, $N(\epsilon,s,n,H)$ and a
collection of balls of radii, $r_i\in [\lambda, \epsilon]$
centered on $z_i\in K$ such that \be \label{packing1}
\bigcup_{i=1}^N B_{z_i}(r_i) \supset K
\textrm{ and }
\sum_{i=1}^N \vol(B_{z_i}(r_i)) \le (1+\epsilon)\vol(K)
\ee
and there are $r_i'<r_i$ such that $B_{z_i}(r_i')$ are disjoint
subsets of $K$ such that
\be \label{packing2}
\sum_{i=1}^N \vol(B_{z_i}(r_i')) \ge (1-\epsilon)\vol(K).
\ee
\end{lem}

\Pf
Fix $\beta<1$.
Let $N_0$ be the maximum number of disjoint balls
contained in $K_0=K$
of radius $\epsilon_0=\epsilon$.
Let
$S_0=\{z_1,...,z_{N_0}\}$
be the centers of these balls and $r_z=\epsilon_0$ their radii.

Since $z_i \in K \subset K_{R_0}$, so $\epsilon < r - R_0 \le r - d_M (p, \pi
 (z_i))$, thus we can apply (\ref{chco1.2}) of
Corollary~\ref{volcor} with an inner ball $B(z_i,\epsilon)$ and an
outer ball up to radius $R-r$.

First we apply it to show that $N_0=N_0(\epsilon,s, n, H)$ does
not depend on the manifold using the standard packing. Recall
that $2s < R-r$.
\be
\vol(B(z_i,\epsilon))\ge
 \frac{V(n,H,\epsilon)}{V(n,H,2s)} \vol(B(z_i,2s))
\ge \frac{V(n,H,\epsilon)}{V(n,H,2s)} \vol(B(\tilde{x},s)).
\ee
So
\be
\vol(B(\tilde{x},s))\ge \sum_{i=1}^{N_0}
\vol(B(z_i,\epsilon)) \ge N_0
\frac{V(n,H,\epsilon)}{V(n,H,2s-\epsilon)}
Vol(B(\tilde{x},s)). \ee

Second we apply (\ref{chco1.2}) to estimate how much of the volume
of $B(\tilde{x},s)$ has been covered by these balls. Note
that $2r_z=2\epsilon < r < R-r$.
 \begin{eqnarray} \label{pack1}
\sum_{i=1}^{N_0} \vol(B_{z_i}(r_z)) & \ge &
\frac{V(H,n,\epsilon)}{V(H,n,2\epsilon)}
\sum_{i=1}^{N_0} \vol(B_{z_i}(2r_z)) \\
& \ge &
\frac{V(H,n,\epsilon)}{V(H,n,2\epsilon)} \vol(K).
\end{eqnarray}

Let
\be
K_{j}=K_{j-1} \setminus \bigcup_{z\in S_{j-1}} B_{z}(r_z)
\ee
Let $N_j$ be the maximum number of disjoint balls
contained in $K_j$ of radius $\epsilon_j=\beta^j\epsilon$.
Let
$S_j=\{z_{N_0+\cdots+N_{j-1}+1}, ..., z_{N_0+\cdots+N_j}\}$
be the centers of these balls and $r_z=\epsilon_j$ be their
radii. As argued above, we can apply (\ref{chco1.2}) with
an inner ball $B(z_i, \epsilon_j)$
because these balls are contained in $B(\tilde{x},s)$ and
apply this to prove $N_j=N_j(\epsilon, \beta, s, n,H)$.

Now, balls of twice the radius cover $K_j$ and we have
\begin{eqnarray} \label{pack2}
\sum_{i=N_0+...+N_{j-1}+1}^{N_0+...+N_j}
\vol(B_{z_i}(r_z))&\ge&
\frac{V(H,n,\beta^j\epsilon)}{V(H,n,2\beta^j\epsilon)}
\sum_{i=N_0+...+N_{j-1}+1}^{N_0+...+N_j}
\vol(B_{z_i}(2r_z))\\ &\ge&
\frac{V(H,n,\beta^j\epsilon)}{V(H,n,2\beta^j\epsilon)}
\vol(K_j) = C(H,n, \beta, \epsilon, j) \vol(K_j),
\end{eqnarray} where
\be \label{CJay}
C(H,n, \beta, \epsilon, j)=
\frac{V(H,n,\beta^j\epsilon)}{V(H,n,2\beta^j\epsilon)}.
\ee

For fixed $\epsilon, \beta<1$ we can take $J_0$ sufficiently large that
$C(H,n, \beta, \epsilon, J)$ is approximately $(1/2)^n$.
More precisely we can take $J_0$ sufficiently large that
\be \label{Jay0}
C(H,n,\beta,\epsilon,j)\ge (1/3)^n \textrm{ for all } j\ge J_0.
\ee
This and the definition of $K_j$ gives us
\begin{eqnarray} \label{pack3}
\vol(K_j) &=& \vol(K_{j-1})-
\sum_{i=N_0+...+N_{j-2}+1}^{N_0+...+N_{j-1}}\vol(B_{z_i}(r_{z
_ i}))  \\ &\le& \vol(K_{j-1})- C(H,n,\beta, \epsilon, j-1)
\vol(K_{j-1}) \\ &\le& \left(1- C(H,n,\beta,\epsilon,j-1)
\right) \vol(K_{j-1}) \\ &\le&
\prod_{k=1}^{j-1}\left(1-C(H,n,\beta,\epsilon,k)\right)
\vol(K) \\ &\le& (1-(1/3)^n)^{(j-J_0)} \vol(K) \qquad \forall
j\ge J_0. \label{3.33} \end{eqnarray}

We must take $J\ge J_0$ sufficiently large that
\be \label{Jay}
3^n(1-(1/3)^n)^{(J-J_0) }      < \epsilon.
\ee

Let $r_i=r_i'=r_{z_i}$ for $z_i\in \bigcup_{j=1}^{J-1} S_j$ and
$r_i=2r_i'=2r_{z_i}=2\beta^J\epsilon$ for $z_i\in S_J$ so that
$$
\bigcup_{i=1}^{N_1+...N_J} B_{z_i}(r_i) \supset K
$$
but the $r_i'$ balls are disjoint, giving us
\begin{equation} \label{Pack4}
\sum_{i=1}^{N_1+...+N_J} \vol(B_{z_i}(r_i')) \le \vol(K).
\end{equation}

Thus, by (\ref{chco1.2}) and the definition of the radii,
the fact that $B_{z_i}(r_i)$ for $i=1, \cdots, N_1+...+N_{J-1}$
 are disjoint and in $K$, and
(\ref{3.33}), we get
\begin{eqnarray} \label{pack5}
\sum_{i=1}^{N_1+...+N_J} \vol(B_{z_i}(r_i)) &\le&
\sum_{i=1}^{N_1+...+N_{J-1}} \vol(B_{z_i}(r_i))+
\sum_{i=N_1+...+N_{J-1}+1}^{N_1+...+N_J} \vol(B_{z_i}(r_i))\\
& \le&
\vol(K) +
\frac{V(n,H,2\beta^J\epsilon)}{V(n,H,\beta^J\epsilon)}
\sum_{i=N_1+...+N_{J-1}+1}^{N_1+...+N_J} \vol(B_{z_i}(r_i'))
\\ &\le&  \vol(K) +
\left(\frac{V(n,H,2\epsilon\beta^J)}{V(n,H,\epsilon\beta^J)}\right) \vol(K_J)\\ &\le&
 \vol(K) + 3^n
(1-(1/3)^n)^{(J-J_0) }\vol(K)
\\
&\le &
 \left(1 + 3^n(1-(1/3)^n)^{(J-J_0)} \right)
\vol(K)   \\
&\le& (1+\epsilon)\vol(K).
\end{eqnarray}

Finally by (\ref{pack3}), we also have
\begin{eqnarray}
\sum_{i=1}^{N_1+...+N_J} \vol(B_{z_i}(r_i'))
& = &
\vol(K)-\vol(K_{J+1}) \\
&\ge&
\vol(K)- (1-(1/3)^n)^{(J+1-J_0)} \vol(K)  \\
&\ge& (1-\epsilon)\vol(K).
\end{eqnarray}

Recall that we can set $\beta=1/2$.
Then $\lambda(\epsilon,s,n,H)=\epsilon\beta^J$
where $J$ is determined in (\ref{Jay}), ({\ref{Jay0}) and (\ref{CJay})
and $N(\epsilon,s,n,H)=\sum_{j=1}^J N_j(\epsilon, \beta, s,n,H)$.
\qed

We can now complete our construction of the measure.
\newline

\noindent {\bf Proof of Proposition~\ref{limmeas}:}
We are given a pointed Gromov Hausdorff converging sequence
of spaces $\tilde{B}(p_i,r_i,R_i)^{\delta_i}$ which
converge as covers to $\hat{B}(y,r,R)$.  We may need to take
a further subsequence to get a renormalized limit measure.

First define the renormalized volume functions
\be
\bar{V}_i:\tilde{B}(p_i,r_i,R_i)^{\delta_i} \times \mathbb R
\to \mathbb R, \qquad \qquad
\bar{V}_i(\tilde{x},\rho)=\vol(B(\tilde{x},\rho))/\vol(\tilde{p}_i,r/10), \ee
where $B(\tilde{x},\rho)$ is defined
using $d_{\tilde{B}(p_i,r_i,R_i)^{\delta_i}}$. As in the
proof of \cite[Theorem 1.6]{ChCo2}, we will show these
$\bar{V}_i$ are uniformly equicontinuous and bounded, but to
do so here we must restrict our domain considerably.

Let $\tilde{p}_x=g\tilde{p}_i$ closest to $\tilde{x}$, then
$d(\tilde{p}_x,\tilde{x})=d(p_i,\pi_i(\tilde{x}))\le r$ and
\be
\bar{V}_i(\tilde{x},\rho)=\vol(B(\tilde{x},\rho))/\vol(\tilde
{ p }_x , r/ 10). \ee

{\em Uniformly Bounding the $\bar{V}_i$:}

Temporarily fix $R_0<r$.  We have $R_0<r_i$ eventually, so we can define
 $K^i_{R_0}=\{\tilde{x}\in \tilde{B}(p_i,r_i,R_i)^{\delta_i}:
d(p_i,\pi_i(\tilde{x})\le R_0\}$ and $K_{R_0}=\{\hat{x}\in
\hat{B}(y,r,R): d(y, \pi(\hat{x}))\le R_0\}$ as in
(\ref{KR0}).  Note that by the given continuity of the
covering maps as $i\to\infty$, we have \be \label{PetKR0}
GH\lim_{i\to\infty} (K^i_{R_0}, \tilde{p}_i) =K_{R_0}.
\ee

For any $R_1\in (0, \min\{ r-R_0, R_0+\frac{r}{10}\})$, and
$R_2\in (R_1, R-r)$, then for $i$ sufficiently large we have
$R_1<r_i-R_0$, $R_1<R_0+r/10<R_i-r_i$, $R_2+R_0 < R_i$.

We restrict $\bar{V}_i(\tilde{x},\rho)$ to $K_{R_0}\times
[R_1,R_2]$. These functions are nondecreasing for fixed
$\tilde{x}$. We can apply this fact and (\ref{chco1.3})
in Corollary~\ref{volcor}, to get a uniform lower bound,
\begin{eqnarray} \bar{V}_i(\tilde{x},\rho) &\ge&
\bar{V}_i(\tilde{x},R_1)  =
\vol(B(\tilde{x},R_1))/\vol(\tilde{p}_x,r/10)\\ &\ge&
V(n,H,R_1)/V(n,H,R_0+r/10)>0, \end{eqnarray}
because  $R_1\le (R_0+r/10) \le R_i-r_i$ and
$R_1<r_i-R_0$ as required by (\ref{chco1.3}).

We get a uniform upper bound for $\bar{V}_i$
on $K_{R_0}\times [R_1,R_2]$
by applying containment and (\ref{volcor1}),
\begin{eqnarray}
\bar{V}_i(\tilde{x},\rho)& \leq
&\bar{V}_i(\tilde{x},R_2) =
\vol(B(\tilde{x},R_2))/\vol(\tilde{p}_x,r/10)\\
&\le&\vol(B(\tilde{p}_x,R_2+R_0))/\vol(\tilde{p}_x,r/10) \\
&\le& V(n,H,R_2+R_0)/V(n,H,r/10). \end{eqnarray}

{\em Equicontinuity of the $\bar{V}_i$:}

We now further restrict the domain to
$K_{R_0}\times [R_1,R_2]$ if $R_2< r-R_0$
so that $\bar{V}_i(x,\rho)$ are uniformly continuous
 in the sense of \cite{GP}, see also \cite[Page 279]{Pe2}.
The restriction on $R_2$ comes from the trouble with estimating
large balls using the intrinsic metric.
Again we take $i$ sufficiently large that $R_2<r_i-R_0$.

Given $\tilde{x}_i \in K^i_{R_0}$ and $\rho_i\in
[R_1,R_2]$. Let
$d_{\tilde{B}(p_i,r,R)^\delta}(\tilde{x}_1,\tilde{x}_2)=s<\delta < r $ . Then $\rho_i<r-d_M(\pi(\tilde{x}_i), p)$ and
$\rho_i < r<R-2r<R-r-d(\tilde{x}_1,\tilde{x}_2)$, so we can
apply (\ref{chco1.2}) to get the following: \ba
\label{diffcor2}
\lefteqn{
|\bar{V}_i(\tilde{x}_2,\rho_2)-\bar{V}_i(\tilde{x}_1,\rho_1)|
=
\frac{|\vol(B(\tilde{x}_2,\rho_2))-\vol(B(\tilde{x}_1,\rho_1)
) | } {\vol B(\tilde{p},r/10)}} \nonumber \\ &=&\frac{\vol
(B(\tilde{x}_2,\rho_2)\setminus B(\tilde{x}_1,\rho_1)) + \vol
(B(\tilde{x}_1,\rho_1)\setminus B(\tilde{x}_2,\rho_2))}{\vol
B(\tilde{p},r/10)} \nonumber \\ & \le &
\frac{\vol(Ann_{\tilde{x}_2}(\rho_2, \rho_1
+s))+\vol(Ann_{\tilde{x}_1}(\rho_1, \rho_2 +s))}{\vol
B(\tilde{p},r/10)} \nonumber \\
&=&\frac{\vol(B(\tilde{x}_2,\rho_2))}{\vol
B(\tilde{p}_{x_2},r/10)}\left(\frac{\vol
(B(\tilde{x}_2,\rho_1 +s))}{\vol(B(\tilde{x}_2,\rho_2))} -1
\right) + \frac{\vol(B(\tilde{x}_1,\rho_1))}{\vol
B(\tilde{p}_{x_1},r/10)}\left(\frac{\vol
(B(\tilde{x}_1,\rho_2+s))}{\vol(B(\tilde{x}_1,\rho_1))} -1
\right) \nonumber \\  & \le & \frac{V(n,H,R_2)}{V(n,H,r/10)}
\cdot
\left(\frac{V(n,H,\rho_1+s)-V(n,H,\rho_2)}{V(n,H,\rho_2)} +
\frac{V(n,H,\rho_2+s)-V(n,H,\rho_1)}{V(n,H,\rho_1)}\right).
\ea Here $Ann_{x}(\rho_1,\rho_2) = B_x (\rho_2) \setminus B_x
(\rho_1)$ and is empty if $\rho_2 < \rho_1$. Since $V(n,H,s)$
is continuous in $s$, we know that for any $\epsilon$ we can
find $\rho_1$ near $\rho_2$ and $s$ small enough that this
last line is less than $\epsilon$.

This gives us uniform equicontinuity with the restricted $R_2<r-R_0$.
Thus we can apply a generalized version of the
Arzela-Ascoli theorem (\cite{GP}, \cite[Page 279, Lemma
1.8]{Pe2}) combined with (\ref{PetKR0}) to get a subsequence
of the $\bar{V}_i$ converging uniformly to a limit function
defined on $K_{R_0} \times [0, r-R_0]$,
\be  \label{vbar}
\bar{V}_\infty(\hat{x},s)=\lim_{i\to\infty}\bar{V}_i(\tilde{x}_i, s). \ee

{\em Extending the domain of $\bar{V}_\infty$:}

Recall the definition of $i\tilde{B}$ in (\ref{KR0}).
Then we have,
\be \label{defnU}
U=\{(\hat{x},s):\hat{x}\in i\hat{B}(y,r,R), s\le
r-d(\hat{x},\hat{p}_x)\} =
\bigcup_{R_0<r} K_{R_0}\times [0, r-R_0].
\ee
We can extend the definition of $\bar{V}_\infty$
by taking a sequence of $R_0\to r$ and diagonalizing the subsequences
used to define $\bar{V}_\infty$ on each $K_{R_0}\times [0,r-R_0]$.

Applying (\ref{chco1.2}) of Corollary~\ref{volcor}, we know
\be
\frac
{\bar{V}_\infty(\hat{x},r_1)}{\bar{V}_\infty(\hat{x},r_2)}
\ge \frac{V(n,H,r_1)}{V(n,H,r_2)} \ \ \forall
r_2<r-d(\hat{x},\hat{p}_x). \ee Note that unlike
(\ref{chco1.2}) we only have this estimate on small balls
because we only had equicontinuity on the small balls'
volumes. Furthermore we have no estimate for volumes of balls
centered on the boundary. This is in strong contrast to
\cite[Theorem 1.6]{ChCo2}.

{\em The Caratheodory Construction:}

We will now construct the renormalized limit measure on all
of $i\hat{B}(y,r,R)$ using a standard Caratheodory
construction as in 2.10 of \cite{Fed} or Method II in
\cite{Mun}.  This is different than the construction used by
Cheeger-Colding because the balls have variable size but
agrees with their construction on sets contained in a
$K_{R_0}$ for $R_0<r$.

We first choose our family of open sets,
\be \label{FBall}
F:=\{B_{\hat{x}}(s): (\hat{x},s)\in U \}
\ee
 where $U$ is defined in (\ref{defnU}).
Note that these balls can be measured using
$d_{\hat{B}(y,r,R)}$ without any difficulties involving the
boundary.

For any $A\subset i\hat{B}(y,r,R)$, let
\be \label{GBall}
G_{\epsilon,A}=\{ S \subset F: \bigcup_{B\in S} B \supset A
\textrm{ and } diam{B}\le\epsilon \,\,\forall B\in S\}.
\ee
The families of open balls in $G_{\epsilon, A}$ may all be
infinite if $A$ is not compact.

Let
\be \label{psiball}
\psi:F \to [0,\infty)\ \mbox{be defined}\
\psi(B(\hat{x},s))=\bar{V}_\infty(\hat{x},s).
\ee

Then define a measure, $\mu$, as  
\be \label{carath2}
\mu(A)=\lim_{\epsilon\to0} \mu_\epsilon(A),
\ee
where
\be \label {carath1}
\mu_\epsilon(A)=\inf_{S\in G_{\epsilon,A}} \{ \sum_{B\in S}
\psi(B)\}. \ee

Since all the members of $F$ are Borel sets, $\mu$ is a Borel
Regular measure \cite[2.10]{Fed}.  This measure is Radon when restricted to
measurable sets by Thm 13.7 of \cite{Mun} because
the members of $F$ are open sets.

{\em Properties of $\mu$: }

If $B(\hat{x},s)\in F$ then
\be \label{limmeas4}
\mu(B(\hat{x},s)) \ge \mu_{2s}(B(\hat{x},s)) =
\psi(B(\hat{x},s))=\bar{V}_\infty(\hat{x},s).
\ee

Now we want to bound $\mu(B(\hat{x},s))$ from above
for {\em arbitrary} $\hat{x}\in i\hat{B}(y,r,R)$
and {\em large} $s< (R-r)/2$ where these balls are measured using
$d_{\hat{B}(y,r,R)}$ not $d_{i\hat{B}(y,r,R)}$.  To avoid
trouble, we use the property of Radon measures \cite[Defn 2.2.5]{Fed},
that \be
\mu(B(\hat{x},s))=\sup\{\mu(\bar{K}),
\bar{K}\subset B(\hat{x},s) \,\,  s.t. \ \bar{K} \textrm{
compact.}\} \ee
Note that for any compact $\bar{K}$ in the open set
$i\hat{B}(y,r,R)$ there exists
$\epsilon=\epsilon_{\bar{K}}>0$ such that
$T_{\epsilon}(\bar{K})$ avoids the boundary of
$i\hat{B}(y,r,R)$. So in fact $\bar{K}\subset
B(\hat{x},s)\cap K_{R_0}$ where $R_0=r-\epsilon_{\bar{K}}$
and $K_{R_0}$ is defined in (\ref{KR0}). Thus in fact, \be
\label{supKR0} \mu(B(\hat{x},s))=\sup\{\mu(B(\hat{x},s)\cap
K_{R_0}), R_0<r\}. \ee

Fix $r/10 <R_0<r$ and set
\be \label{KBxs}
K=B(\hat{x},s)\cap K_{R_0}.
\ee
Let $\tilde{B}(p_i,r_i,R_i)^{\delta_i}$ be within
$\epsilon_i<\epsilon_K/10=(r-R_0)/10$ of $\hat{B}(y,r,R)$,
and let \be \label{KiB}
K_i=B(\tilde{x}_i,s+\epsilon_i)\cap
K^i_{R_0+\epsilon_i}\subset
\tilde{B}(p_i,r_i,R_i)^{\delta_i}. \ee
Now applying (\ref{volcor1}) as in (\ref{diffcor2}), and noting
that $s+\epsilon_i \le \frac{R_i-r_i}{2}$ for all $i$ large, we have
\begin{eqnarray}
\frac{|\vol(K_i) - \vol(B(\tilde{x}_i,s+\epsilon_i))|}
{\vol(B(\tilde{p}_i, r/10))} &\le & \sum_{g\tilde{p}_i \in
B(\tilde{x}_i, s+\epsilon_i)}
\vol(Ann_{g\tilde{p}_i}(R_0,r_i))/\vol(B(g\tilde{p}_i, r/10))
\nonumber \\
& = & \sum_{g\tilde{p}_i \in B(\tilde{x}_i, s+\epsilon_i)}
\frac{\vol(B(\tilde{p}_i, R_0))}{\vol(B(\tilde{p}_i,
r/10))}\left( \frac{\vol(B(\tilde{p}_i,
r_i))}{\vol(B(\tilde{p}_i, R_0))} -1\right) \nonumber \\
 &\le &  \#\{g\tilde{p}_i \in B(\tilde{x}_i, s+\epsilon_i)\}
\frac{V(n,H,r)-V(n,H,R_0)}{V(n,H,r/10)}
\end{eqnarray}
for all $i$ large. By the given isometry of $\delta/2$
balls of the covering maps and (\ref{chco1.3}),
$\#\{g\tilde{p}_i \in B(\tilde{x}_i, s+\epsilon_i)\} \leq
\frac{V(n,H,s+\epsilon_i)}{V(n,H,\delta/2)}\leq
\frac{V(n,H,2s)}{V(n,H,\delta/2)}$.

So the left side of
this equation is uniformly bounded for all $i$ large, therefore
\be \label{VolKi}
\forall \varepsilon>0, \,\, \exists R_0 \textrm{
sufficiently close to } r, \textrm{ s.t. }
\left|\frac{\vol(K_i)}{\vol(B(\tilde{p}_i,
r/10))}-\bar{V}_i(\tilde{x}_i,s+\epsilon_i)\right|\le
\varepsilon \frac{V(n,H,2s)}{V(n,H,\delta/2)}.
\ee
for all $i$ sufficiently large.

By Lemma~\ref{packing}, for
all $\epsilon>0$, there exists $\lambda$ and $N$ depending
only on $\epsilon$, $s+\epsilon_K$, $n$ and $H$ such that $K_i$ has two
special families of balls, $H_i$ and $\bar{H}_i$ consisting
of at most $N$ balls each.  These balls have of radii
between $\epsilon$ and $\lambda$ and satisfy:
\be
T_\epsilon(K_i) \supset
\bigcup_{B\in H_i} B \supset K_i
\textrm{ and }
\sum_{B\in H_i} \vol(B) \le (1+\epsilon)\vol(K_i)
\ee
and $B$ in $\bar{H}_i$ are disjoint subsets of $K_i$
such that
\be
\sum_{B\in \bar{H}_i} \vol(B) \ge (1-\epsilon)\vol(K_i).
\ee

Taking $\epsilon<\epsilon_K/10$ guarantees that if
$B(\tilde{x}',s')\in H_i$ then there exists $R_0<r$ such
that $(\tilde{x}',s')\in K_{R_0}^i \times [0, r-R_0]$.
Thus
\be
\vol(B(\tilde{x}',s'))=\bar{V}_i(\tilde{x}',s')
\vol(B(p_i,r/10)) \ \ \forall (\tilde{x}',s')\in U. \ee

Let $\Phi_i:\tilde{B}(p_i,r_i,R_i)^{\delta_i} \to
\hat{B}(y,r,R)$ be the $\epsilon_i$ Hausdorff approximation.

Let $F_i'=\{B_{\Phi_i(\tilde{x}')}(s'+2\epsilon_i): \,\,
B(\tilde{x}',s')\subset H_i\}$. By the choice of $\epsilon$
and $\epsilon_i$, we know $F_i' \subset F$.
In fact if $B(\hat{x}',s')\in F_i'$ then
$(\hat{x}',s')\in K_{r-\epsilon}\times [\lambda,\epsilon]$.
Furthermore for each $i$ sufficiently large,
\be
K \subset \bigcup_{B\in F_i'} B.
\ee
Thus by (\ref{carath1}) and the uniform convergence of $\bar{V}_i$
there exists $\bar{\epsilon}_i \to 0$,
\begin{eqnarray}
\mu_{\epsilon+2\epsilon_i}(K)& \le& \sum_{B(\hat{x}',s')\in
F_i'} \bar{V}_\infty(\hat{x}',s') \\
& \le & \sum_{B(\tilde{x}',s')\in H_i}
\bar{V}_i(\tilde{x}',s') + N \bar{\epsilon_i} \\
& = & \sum_{B(\tilde{x}',s')\in
H_i} \vol(B(\tilde{x}',s'))/\vol(B(\tilde{p}_i,r/10)) + N
\bar{\epsilon_i}\\ &\le &
(1+\epsilon)\vol(K_i)/\vol(B(\tilde{p}_i,r/10)) + N
\bar{\epsilon_i}. \end{eqnarray} Applying our estimate on the
volume of $K_i$ in (\ref{VolKi}), we get \be \label{limmeas3}
\mu_{\epsilon+2\epsilon_i}(K) \le
 (1+\epsilon)\left(\bar{V}_i(\tilde{x}_i,s+\epsilon_i)+
\varepsilon \frac{V(n,H,2s)}{V(n,H,\delta/2)}\right) +N
\bar{\epsilon_i}. \ee
We will apply this equation for large and small $s$.

First we look at small $s<r-d(\hat{x},g\hat{p})$.  Since
$N$ depends on $\epsilon$ but not on $i$, taking $i$ to
infinity, $\epsilon_i,\bar{\epsilon}_i$ go to $0$, and we get
\be
\mu_{\epsilon}(K)\le
(1+\epsilon)\left(\bar{V}_\infty(\hat{x},s)+
\varepsilon \frac{V(n,H,2s)}{V(n,H,\delta/2)}\right).
\ee
Then taking $\epsilon$ to zero, we get
\be
\mu(K)\le \bar{V}_\infty(\hat{x},s)+\varepsilon
\frac{V(n,H,2s)}{V(n,H,\delta/2)}. \ee
Taking the supremum over $K=B(\hat{x},s)\cap K_{R_0}$, using
(\ref{supKR0}) with $R_0 \to r$
and (\ref{KBxs}), and then taking $\varepsilon$ to zero
as in (\ref{VolKi}) we get:
\be
\mu(B(\hat{x},s))\le \bar{V}_\infty(\hat{x},s).
\ee
Combining this with (\ref{limmeas4}), we get
\be \label{limmeas5}
\mu(B(\hat{x},s))=\bar{V}_\infty(\hat{x},s), \qquad \forall
B(\hat{x},s)\in i\hat{B}(y,r,R),
\ee
which gives us (\ref{vinfeq}).

To examine $s<(R-r)/2$, we use (\ref{limmeas3}) again.
We apply the volume comparison for $s$, to get
for $r_1<r_i-d(\tilde{x}_i,g\tilde{p_i})$, $r_1<s$,
\be
\mu_{\epsilon+2\epsilon_i}(K) \le (1+\epsilon)
\left(\frac{V(n,H,s+\epsilon_i)}{V(n,H,r_1)}\bar{V}_i(\tilde{x}_i, r_1) + \varepsilon
\frac{V(n,H,2s)}{V(n,H,\delta/2)} \right)
+N \bar{\epsilon_i}. \ee

Now taking $i$ to infinity and $\epsilon_i,\bar{\epsilon}_i$
to 0, and lastly $\epsilon$ to zero, we get
\be
\mu(K)\le
\frac{V(n,H,s)}{V(n,H,r_1)}\bar{V}_\infty(\hat{x},r_1)
+ \varepsilon \frac{V(n,H,2s)}{V(n,H,\delta/2)}. \ee
Taking the supremum over $K=K_{R_0}\cap B(\hat{x},s)$ with
$R_0 \to r$ as in (\ref{supKR0}) and (\ref{KBxs}), and
finally taking $\varepsilon$ to zero as in (\ref{VolKi}) we
get: \be \label{limmeas6}
\mu(B(\hat{x},s))\le
\frac{V(n,H,s)}{V(n,H,r_1)}\mu(B(\hat{x},r_1)) \ee
which gives us (\ref{mu1.2}).

We now define $\mu$ as a measure on $\hat{B}(y,r,R)$, by
setting $\mu(A)=\mu(A\cap i\hat{B}(y,r,R))$ for any Borel set
$A$. Then (\ref{limmeas5}) implies (\ref{vinfeq})
and (\ref{limmeas6}) implies (\ref{mu1.2}) for this $\mu$
and $\mu$ is still Borel regular.
\qed

\subsection{Stability}

We can now use the limit covers $B(y,r,R)^\delta$ and
their measures [Theorems~\ref{limdelcov} and~\ref{limmeas}]
to get the stability required by Theorem~\ref{univstab} to
prove the existence of a universal cover for $Y$
[Theorem~\ref{mainthm}].

\begin{theo}  \label{samedelta}
For all $R>r>0$, $y \in Y$ there exists $\delta_{y,r,R}$ depending on
$Y, y, r,R$ such that for all $\delta<\delta_{y,r,R}$, we have
\be \label{sdcoverproj}
\tilde{B}(y,r,R)^\delta
={B}(y,r,R)^{\delta}=\tilde{B}(y,r,R)^{\delta_{y,r,R}}.
\ee
\end{theo}

Note there is no restriction on $R$ and $r$ in this statement.

To prove this we will first prove that for special regular points $y\in Y$
(which are proven to be dense in $Y$ by Cheeger-Colding),
sufficiently small balls lift isometrically to all covers .
Recall that a regular point is a point in a metric space whose
tangent cone is Euclidean and therefore has a pole.

\begin{theo}  \label{regdel}
Let $(Y,p)$ be the pointed Gromov-Hausdorff limit of a sequence of complete
manifolds $(M^n_i, p_i)$
such that
\be \label{cond3}
\Ric_{M_i} \ge -(n-1)H \textrm{ where } n\ge 3, H>0.
\ee
If $y \in B(p,1) \subset Y$ is a point such that
there exists a tangent cone,
($Y^\infty, y_\infty$),
that has a pole at $y_\infty$, then for any
$100 \le 10 \bar r \le \bar R$ there
exists $r_y(\bar r,\bar R)>0$, such that for
all $\delta>0$, $B(y,r_y)$ lifts isometrically
to $B(p,\bar r,\bar R)^\delta$.
\end{theo}

The proof of Theorem~\ref{regdel}
uses the Abresch-Gromoll Excess estimate on the relative
$\delta$ cover as in \cite[Thm 4.5]{SoWei}
except that now our covers have boundary.
In particular, we need the following adaption of the
Abresch Gromoll Excess Theorem \cite{AbGl} for manifolds
with boundary.

\begin{lem}  \label{ptexcess}
Let $M^n$ be a compact Riemannian manifold with boundary satisfying
(\ref{cond3}). For a ball $B(p,10\rho) \subset M$ not intersecting
$\partial M$,
there exists a constant
\be \label{EqnSnk}
S=S_{n, H}=\min \left\{  \frac{1}{8},
\frac{1}{4\cdot 3^{n}} \frac{1}{\cosh(\sqrt{H}/4)} \frac{n}{n-1}
\left(\frac{n-2}{n-1} \right)^{n-1}
\left( \frac{\sqrt{H}}{\sinh \sqrt{H}} \right)^{n-1} \right\}
\ee
such that if  $\gamma \subset B(p,\rho)$
is a length minimizing curve of $M$ with
length $D\le 1$
and $x \in B(p,\rho)$ satisfying
$$
d_M(x, \gamma(0))\ge (S_{n,H}+1/2)D \ \textrm{ and } \
d_M(x, \gamma(D))\ge (S_{n,H}+1/2)D,
$$
then
\be \label{abgl0}
d_M(x, \gamma(D/2)) \ge 3S_{n,H}D.
\ee
\end{lem}

This lemma holds as in the proof of Lemma 4.6 in \cite{SoWei} except
for a small technicality involving the use of intrinsic versus
restricted metrics. To overcome this, one notes that
$d_M(x, \gamma (0)) \leq 2\rho$, so
$$
B_{d_M}(x, 2\rho) \subset B_{d_M}(\gamma (0), 4\rho) \subset
B_{d_M}(p, 5\rho).
$$
By Lemma~\ref{distance}, $d_{M}$ and $d_{B(p,10\rho)}$
restricted to $B(p, 5\rho)$ are same.

\vspace{.2in}
\noindent{\bf Proof of Theorem~\ref{regdel}}.
Since $y$ is close enough to $p$, all the points and curves involved in the
proof of \cite[Theorem 4.5]{SoWei}
lies in $B(p,4.8)$, far away from the boundary
$\partial B(p, \bar r)$, similarly for the cover. By
Lemma~\ref{distance}, the restricted distance on $B(p,4.8)$ from
$B(p,\bar r)$ and from $Y$ are same. So the proof of
\cite[Theorem 4.5]{SoWei} carry over.




Assume on the contrary that for all $r>0$ there is a $\delta_r>0$
such that the ball $B(y,r)$ does not lift isometrically to $B(p,\bar r,\bar
R)^{\delta_r}$. Let $G^\delta$ denote the deck transformation group on $B(p,\bar
r,\bar R)^\delta$. Thus,
there exist $r_i \to 0$, $\delta_i=\delta_{r_i}$,
and $g_i\in G^{\delta_i}$ such that $d_i=d_{B(p,\bar r,\bar
R)^{\delta_i}}(\tilde{y}, g_i\tilde{y}) \in (0, 2r_i) \subset
(0,1]$. In fact, we can choose $g_i$ so that
\be \label{shortest}
d_{B(p,\bar r,\bar R)^{\delta_i}}(\tilde{y}, g_i\tilde{y})
\le d_{B(p,\bar r,\bar R)^{\delta_i}}(\tilde{y}, h\tilde{y}) \ \
\forall h\in  G^{\delta_i}.
\ee

Next we will find a length minimizing curve, $\tilde{C}_i$, running from
$\tilde{y}$ to  $g_i\tilde{y}$ which has the property that it
passes through a particular point $\tilde{z}_i= \tilde{C}_i(d_i/2)$
which is the limit of halfway points of length minimizing curves in
the sequence $\tilde{B}(p_i,\bar r_i, \bar R_i)^{\delta_i}$.
We do this so that
we can apply
Lemma~\ref{ptexcess} to $\tilde{B}(p_i,\bar r_i, \bar R_i)^{\delta_i}$.

To construct $\tilde{C}_i$, we first let $\tilde{y}_j$, $\tilde{y}^i_j \in
B(\tilde{p_j}, 2.4) \subset \tilde{B}(p_j,\bar r_j, \bar R_j)^{\delta_i}$ which
are close to $\tilde{y}$ and
$g_i\tilde{y}$.  So $d_{\tilde{B}(p_j,\bar r_j, \bar R_j)^{\delta_i}}
(\tilde{y}_j, \tilde{y}^i_j)=d_{i,j}$
converges to $d_i$.  Let $\tilde{z}^i_j$ be midpoints of minimal geodesics
$\gamma^i_j$,
running from $\tilde{y}_j$ to $\tilde{y}^i_j$.  Taking a subsequence
of $j\to\infty$, there is a point $\tilde{z}_i\in B(p,\bar r, \bar
R)^{\delta_i}$ which
is halfway between $\tilde{y}$ to  $g_i\tilde{y}$.  Let $\tilde{C}_i$
be a length minimizing curve running from $\tilde{y}$ to $\tilde{z}_i$
and then to $g_i\tilde{y}$.  Finally let $C_i$ be the projection of
$\tilde{C_i}$ to $B(p,\bar r)$. $C_i$ lies in $B(p,4.8)$.

Now, imitating the proof of the Halfway Lemma of [So],
and using (\ref{shortest}), we
know $C_i\in B(p,4.8)$ is minimizing halfway around, $d_{Y}(C_i(0),
C_i(d_i/2))=d_i/2$.

We choose a subsequence of these $i$ such that $(Y,y)$ rescaled
by $d_i$ converges to a tangent cone $(Y^\infty, y_\infty)$.  So
\begin{equation}  \label{tancone}
d_{GH}\left( B(y,10d_i)\subset Y, B(y_\infty, 10d_i)\right) < \epsilon_i d_i
\end{equation}
where $\epsilon_i$ converges to $0$.

Let $S$ be the constant from Lemma~\ref{ptexcess}.
Since $Y^\infty$ has a pole at $y_\infty$, we know there is a length
minimizing curve running from $y_\infty$ through any point in
$\partial B(y_\infty, d_i/2)$ to $\partial B(y_\infty, d_i/2+2S d_i)$.
Thus by (\ref{tancone}),
\be \label{wherex}
\forall \,\,\,x\in \partial B(y,d_i/2+2S d_i) \subset Y,
\ee
we have points
\be
x_\infty\in Ann_{y_\infty}(d_i/2+2S d_i-\epsilon_id_i,
d_i/2+2S d_i+\epsilon_id_i)
\ee
and
\be
y_i\in Ann_{y_\infty}(d_i/2-\epsilon_id_i,
d_i/2+\epsilon_id_i)
\ee
such that
\begin{eqnarray} \label{epsicont}
d_Y(x, C_i(d_i/2)) &<& d_{Y^\infty}(x_\infty, y_i) + \epsilon_i d_i\\
&\le & 2\epsilon_i d_i +2S d_i +\epsilon_i d_i.
\end{eqnarray}

Now we will imitate the Uniform Cut Lemma of [So], to show that
for all $x\in \partial B(y, d_i/2+2S d_i)$, we have
$l_i=d_Y(x, C_i(d_i/2)) \ge (3S)d_i.$
This will provide a
contradiction for $\epsilon_i < S/2$ and we are done.

First we lift our points $x$ and $y$ to the cover $B(p, \bar r, \bar
R)^{\delta_i}$
as follows.  We lift $y$ to the point $\tilde{y}$ and we lift the closed
loop $C_i$ to the curve $\tilde{C}_i$ running from $\tilde{y}$
through $z_i=\tilde{C_i(d_i/2)}$ to $g_i\tilde{y}$.  Then if $\sigma$
is a length minimizing curve of length $l_i$ running from
$C_i(d_i/2)$ to $x$, we lift it to $\tilde{Y}^{\delta_i}$ so it
runs from $\tilde{z}_i$ to a new point, $\tilde{x}$.
Note that by our choice of $x$ in (\ref{wherex}),
\be
d_{\tilde{Y}^{\delta_i}}(g_i\tilde{y}, \tilde{x})\ge
d_{{Y}}({y}, {x})=d_i/2+2S d_i
\ee
and so is $d_{\tilde{Y}^{\delta_i}}(\tilde{y}, \tilde{x})$.

By our choice of $\tilde{C}_i$ and $\tilde{z}_i$, we know there are
corresponding
points in $\tilde{B}(p_j, \bar r_i, \bar{R}_j)^{\delta_i}$.
That is there is a
triangle
formed by $\tilde{y}_j$, $\tilde{y}^i_j$,
with a minimal geodesic $\gamma^i_j$ running between them  and some point
$\tilde{x}_j$ such that
\begin{eqnarray*}
d_{i,j}=d_{\tilde{B}(p_j, \bar r_j, \bar R_j)^{\delta_i}}(\tilde{y}_j,
\tilde{y}^i_j) & \to& d_i ,\\
d_{\tilde{B}(p_j, \bar r_j, \bar R_j)^{\delta_i}}
(\tilde{y}_j, \tilde{x}_j)&\to&
d_{B(p, \bar r, \bar R)^{\delta_i}}(\tilde{y}, \tilde{x})=(1/2+2S)d_i\\
d_{\tilde{B}(p_j, \bar r_j, \bar R_j)^{\delta_i}}(\tilde{y}^i_j,
\tilde{x}_j)&\to&
d_{B(p, \bar r, \bar R)^{\delta_i}}(g_i\tilde{y}, \tilde{x})=(1/2+2S)d_i.\\
l_{i,j}=d_{\tilde{B}(p_j, \bar r_j, \bar
R_j)^{\delta_i}}(\tilde{\gamma}^i_j(d_{i,j}/2), \tilde{x}_j)
&\to& d_{B(p, \bar r, \bar R)^{\delta_i}}(\tilde{z}_i, \tilde{x})=l_i.
\end{eqnarray*}
So for $j$ sufficiently large, we have
\be
d_{\tilde{B}(p_j, \bar r_j, \bar R_j)^{\delta_i}}(\tilde{y}_j, \tilde{x}_j)\ge
(1/2+S)d_{i,j}
\textrm{ and }
d_{\tilde{B}(p_j, \bar r_j, \bar R_j)^{\delta_i}}(\tilde{y}^i_j, \tilde{x}_j)\ge
(1/2+S)d_{i,j}
\ee
and can apply Lemma~\ref{ptexcess} to get
\be
l_{i,j} \ge 3Sd_{i,j}.
\ee

Taking $j$ to infinity, we get the limit of this bound in
$B(p, \bar r, \bar R)^{\delta_i}$,  namely $l_i \ge 3S d_i$.
This contradicts (\ref{epsicont})
for $\epsilon_i < S/2$ and we are done.
\qed

We can now prove our stability theorem.  We first state a more
geometrically intuitive theorem and then prove that it implies
Theorem~\ref{samedelta}.

\begin{theo} \label{samedelta2}
For all $R>0$, $y \in Y$ one of the following two statements
holds:

I: There exists $\delta_{y,R}$ depending on
$Y, y, R$ such that for all $\delta<\delta_{y,R}$, we have
\be
\tilde{B}(y,R)^\delta=\tilde{B}(y,R)^{\delta_{y,R}}.
\ee

II. For all $R'<R$ there exists $\delta_{R'}$ depending
on $Y,y,R,R'$ such that
\be
\tilde{B}(y,R',R)^{\delta_{R'}}=\tilde{B}(y,R',R)^\delta
\qquad \forall \delta<\delta_{R'}.
\ee
\end{theo}

Note that Theorem~\ref{samedelta2} is essentially saying
that if there is a problem with arbitrarily small noncontractible
curves, then they are near the boundary.

\noindent {\bf Proof of Theorem~\ref{samedelta}:}
For a given $y,r$ and $R$, we apply Theorem~\ref{samedelta2}
to $B(y,R)$.  If case I holds, then taking $\delta_{y,r,R}=\delta_{y,R}$
we have $B(y,r,R)^\delta=B(y,r,R)^{\delta_{y,r,R}}$ for any $r<R$.
If case II holds, we let $R'=r$ and $\delta_{y,r,R}=\delta_{R'}$.
Then (\ref{sdcoverproj}) follows from (\ref{coverproj}) of
Theorem~\ref{limdelcov}.
\qed

\noindent {\bf Proof of Theorem~\ref{samedelta2}:}
Suppose neither I nor II hold.
Then there exists $y$ and $R$ and $\delta_i$ converging to $0$ such
that $\tilde{B}(y,R)^{\delta_i}$ are all distinct.  Then
there exists a sequence of $\delta_i>0$ with
$\delta_1 \leq R/10, \delta_i > 10 \delta_{i+1}$
such that all $\tilde{B}(y,R)^{\delta_i}$  and $G(y,R,\delta_i)$
are distinct.   In particular there are nontrivial elements of
$G(y,R, \delta_i)$ which are trivial in $G(y,R, \delta_{i-1})$.
So there exist $x_i\in B(y,R)$, such that the $B_{x_i}(\delta_{i-1})$ contains
a noncontractible loop, $C_i$, which lifts non-trivially in
$\tilde B(y,R)^{\delta_i}$.

In fact we can choose $x_1$ to be the point
closest to $y$ such that $B_{x_1}(R/10)$ contains a noncontractiblle
loop and then choose $\delta_1\in (0,R/10]$ as small as possible
such that $B_{x_1}(R/10)$ contains a loop $C_1$ which lift nontrivially to
$\tilde B(y,R)^{\delta_1}$.  We can then choose iteratively
$x_j$ the point closest to $y$ such that $B_{x_j}(\delta_{j-1}/10)$
contains a noncontractible loop.  Then set $\delta_{j}\in (0,\delta_{j-1}/10]$
as small as possible so that $B_{x_j}(\delta_{j})$ contains a
loop $C_j$ which lift nontrivially to $\tilde B(y,R)^{\delta_j}$.  Note that
$d_{B(y,R)}(y,x_j)$ is a nondecreasing sequence.

By compactness,
a subsequence of the $x_i$ converge to some point $x$ in $B(y,R)$.

If $x\in \partial B(y,R)$, then for any $R'<R$, we know that
there exists $N_1$ sufficiently large such that
$d_{B(y,R)}(y,x_j)>R'+(R-R')/2$ for all $j\ge N_1$.  There
exists $N_2>N_1$ such that $\delta_{(j-1)/10} <(R-R')/2$ for
all $j\ge N_2$.

By the choice of our sequence of $x_j$,
this implies that if $C$ is a loop contained in $B(y',\delta)$
where $\delta\le \delta_{(N_2-1)/10}$ and
$B(y',\delta)\cap B(y,R')$ is nonempty, then $C$ is contractible
in $B(y,R)$.
Thus
\be
\tilde{B}(y,R',R)^{\delta_{R'}}=\tilde{B}(y,R',R)^\delta  \ \
\forall \delta< \delta_{(N_2-1)/10} = \delta_{R'}.
\ee
This implies Case II which we have assumed to be false.

So now we know $x$ is not in the boundary of $B(y,R)$, and we proceed
to find a contradiction.

Let
$\bar{R}>0$ be defined such that $B(x, \bar{R})\subset B(y,R)$
and let $\bar{r}=\bar{R}/10>0$.  Eventually the $C_i$ are
in $B(x, \bar{r}/6)$.

Note $C_i\in B_{x_i}(\delta_{i-1})$ so they lift as closed curves to
$\tilde{B}(x,\bar{r},\bar{R})^{\delta_{i-1}}$.  Since they lift nontrivially to
$\tilde B(y,R)^{\delta_{i}}$, $C_i$ also lift nontrivially to
$\tilde{B}(x,\bar{r},\bar{R})^{\delta_{i}}$.

Since $C_i$ must lift to a union of balls $B_{g\tilde{x}_i}(\delta_{i-1})$ in
$\tilde{B}(x,\bar{r},\bar{R})^{\delta_i}$, there exists $g_i$ nontrivial in
$G(x,\bar{r},\bar{R},\delta_i)$ such that
\be
d_{\tilde{B}^\delta_i}(g_i\tilde{x}_i, \tilde{x}_i) < 2 \delta_{i-1}.
\ee
Let $\alpha_i$ be the projection of the minimal curve from
$g_i\tilde{x}_i$ to $\tilde{x}_i$.  Then $L(\alpha_i)<2\delta_{i-1}<\bar{r}$
and
\be
\alpha_i\subset B_{x_i}(2\delta_{i-1})\subset B(x, \bar{r}/6+2\delta_i)
\subset B(x, \bar{r}/3).
\ee
The $\alpha_i$ represents an element $g_i$
of ${\pi}_1(B(x,\bar{r}))$ which is mapped non-trivially into
$G(x,\bar{r},\bar{R},\delta_i)$ and
trivially into $G(x,\bar{r},\bar{R},2\delta_{i-1})$.

For any $j$, the limit cover $B(x,\bar{r},\bar{R})^{\delta_j}$ covers
$B(x,\bar{r},\bar{R})^{\delta_i}$ for $i=1... j-1$.  By
Theorem~\ref{limdelcov} $g_1,...,g_{j-1}$ are distinct
nontrivial deck transforms of
$B(x,\bar{r},\bar{R})^{\delta_j}$.

Furthermore, for any $q\in B(x,\bar{r})$,
letting $\tilde{x}_i$ be the lift of $x_i$ closest to
 $\tilde{q}\in B(x,\bar{r},\bar{R})^{\delta_j}$,
we have,
\begin{eqnarray*}
d_{B(x,\bar{r},\bar{R})^{\delta_j}}(\tilde{q}, g_i\tilde{q})
& \le & d_{B(x,\bar{r},\bar{R})^{\delta_j}}(\tilde{q},
\tilde{x}_i)
+d_{B(x,\bar{r},\bar{R})^{\delta_j}}(\tilde{x}_i,
g_i\tilde{x}_i)
+d_{B(x,\bar{r},\bar{R})^{\delta_j}}(g_i\tilde{x}_i,
g_i\tilde{q}) \\ &\le& 2\bar{r} + 2 L(\alpha_i) +2\bar{r} \le
6\bar{r}. \end{eqnarray*}

So the $\delta_j$-length of $g_i$,
defined in Definition~\ref{defnlength} is
\be
l(g_i,\delta_j)=\inf_{q\in B(x,\bar{r})}
d_{B(x,\bar{r},\bar{R})^{\delta_j}}(\tilde{q}, g_i\tilde{q})\le 6\bar{r}.
\ee

Therefore we have for any $j$,
there are $j-1$  distinct elements in
$G(x,\bar{r},\bar{R})^{\delta_j}$ with $l(g_i, \delta_j) \leq 6\bar{r}$.

On the other hand we claim that the total number of elements in
$G(x, \bar{r}, \bar{R})^\delta$ of
$\delta$-length $\leq 6\bar{r}$ is uniformly bounded for all $\delta$ in
terms of geometry and topology of $B(x,\bar{r})$.

To show this claim, let us look
at the lift of a regular point $p \in B(x,\bar{r})$ in the cover
$B(x,\bar{r},\bar{R})^{\delta/2}$.

We know by Theorem~\ref{regdel},
there is a $\delta_p>0$ such that
the ball of radius $\delta_p$ about $p$ is isometrically lifted
to disjoint balls of radius $\delta_p$ in $B(x,\bar{r},\bar{R})^\delta$.
Let
\be
 \delta_0=\min\{\delta_p, \bar{r}\}.
\ee

Let $N$ be the number of distinct elements in
$G(x, \bar{r}, \bar{R})^\delta$ of
$\delta$-length $\leq 6\bar{r}$. Note that $g B(\tilde{p}, \delta_0)$
is contained in
$B(\tilde{p}, 6\bar{r}+\delta_0 )\subset
\tilde{B}(x,\bar{r},\bar{R})^\delta$ for all
$g \in G(x, \bar{r}, \bar{R})^\delta$
with $l(g, \delta) \leq 6\bar{r}$.

Now we cannot control $\tilde{B}(x,\bar{r},\bar{R})^\delta$
very well because it does not have a Bishop Gromov volume
comparison theorem.  However in Theorem~\ref{limmeas} we showed
there is a good measure, $\mu$, on
the limit cover $B(x,\bar{r},\bar{R})^{\delta/2}$, which is a
cover of $\tilde{B}(x,\bar{r},\bar{R})^\delta$ by
Theorem~\ref{limdelcov}. So we can lift the distinct $g$ by
lifting representative minimal curves from $\tilde{p}$ to
$g\tilde{p}$. Thus there are $N+1$ isometric disjoint balls of
radius $\delta_0$ contained in a ball of radius
$6\bar{r}+\delta_0$ in the limit cover.  Here we have included the
center ball as well.

Thus,
applying the properties of $\mu$ from Prop~\ref{limmeas},
we have
\be
N+1 \leq  \frac{\mu(B(\tilde{q}, 6\bar{r}+\delta_0))}
{\mu(B(\tilde{q}, \delta_0))}
\le V(n,H,  6\bar{r}+\delta_0)/V(n,H,\delta_0)
\ee
because $6\bar{r}+\delta_0<\bar{R}-\bar{r}$.

This gives us a contradiction.
\qed

We now prove Corollary~\ref{TKi} concerning the convergence
of measures of sets.

\noindent{\bf Proof of Corollary~\ref{TKi}:}
We  must show given any $\epsilon_i \to 0$,
\be \label{TKi1}
\lim_{i\to\infty}
\frac{Vol(T_{\epsilon_i}(K_i))}{Vol(B(\tilde{p}_i,r/10))} =
\mu(K). \ee
Since $K\subset int(X)$ is compact, there is an $\epsilon_0>0$
such that $T_{4\epsilon_0}(K)\subset int(X)$ and so we can use the
Cheeger-Colding Caratheodory method to define the measure.

Given any $V>0$ and $\epsilon<\epsilon_0$,
there exist $z_1,...,z_n \in K$, such that
$\bigcup_j B_{z_j}(r_j)\supset K$, $r_j<\epsilon$ and
\be
\mu(K)\le \sum_{j=1}^n \mu(B_{z_j}(r_j))\le \mu(K)+V.
\ee
Thus by the definition of the renormalized limit and the
uniform continuity of $\bar{V}_i$,
there are $z_{j,i}\in K_i$, such that
\be
\mu(K)\le \lim_{i\to\infty}\sum_{j=1}^n
\frac{Vol(B_{z_{ji}}(r_j+2\epsilon_i'+\epsilon_i))}
{Vol(B(\tilde{p}_i,r/10))}\le \mu(K)+V
\ee
where $\epsilon_i'$
is the Gromov Hausdorff estimate from $X_i$ to $X$
and $\epsilon_i$ is any sequence converging to $0$.

In particular, for all $i\ge N_1$,
\be
\mu(K)- V \le \sum_{j=1}^n
\frac{Vol(B_{z_{ji}}(r_j))}{Vol(B(\tilde{p}_i,r/10))}\le
 \sum_{j=1}^n
\frac{Vol(B_{z_{ji}}(r_j+2\epsilon_i'+\epsilon_i))}
{Vol(B(\tilde{p}_i,r/10))}\le \mu(K)+2V.
\ee
By the choice of $\epsilon_i'$,
\be
\sum_{j=1}^n Vol(B_{z_{ji}}(r_j+2\epsilon_i'+\epsilon_i))
\ge Vol(T_{\epsilon_i}(K_i)).
\ee
Thus
\be \label{epi}
\limsup_{i\to\infty}\frac{ Vol(T_{\epsilon_i}(K_i))}{Vol(B(\tilde{p}_i,r/10))}
\le \mu(K)+2V \qquad \forall V>0,
\ee
so the limsup is $\le \mu(K)$.
On the other hand
\be
\sum_{j=1}^n Vol(B_{z_{ji}}(r_j))
\le Vol(T_\epsilon(K_i)).
\ee
So for all $\epsilon>0$
\be
\liminf_{i\to\infty}\frac{ Vol(T_{\epsilon}(K_i))}{Vol(B(\tilde{p}_i,r/10))}
\ge \mu(K).
\ee
So we can choose $\epsilon_i \to 0$ such that
\be
\liminf_{i\to\infty}\frac{ Vol(T_{\epsilon_i}(K_i))}{Vol(B(\tilde{p}_i,r/10))}
\ge \mu(K).
\ee
and using the same $\epsilon_i$ in (\ref{epi}), we are done.
\qed

\sect{Properties of $\tilde{Y}$ and Applications}\label{appl}

In this section we first study properties of the universal
cover $\tilde{Y}$.  We begin by showing that domains in the
universal cover are Gromov Hausdorff limits of relative delta covers
[Theorem~\ref{univisalim}]
and [Corollary~\ref{2univisalim}].  We then easily show that there is a
measure on $\tilde{Y}$ such that Bishop-Gromov's volume comparison holds
[Theorem~\ref{VOL}].
We then show that the Cheeger-Gromoll Splitting Theorem [Theorem~\ref{SPLIT}]
holds.
Finally we will derive some applications.

Recall, as in Example~\ref{cylinder},
that there may be no covers of the $M_i$ that
converge to $\tilde{Y}$.  The following example demonstrates that a
region in $\tilde{Y}$
may not be limit of relative delta covers with a uniformly
bounded outer radius, $R$.

\begin{ex} \label{cooper}
There is a complete noncompact
$4$-manifold $M$ with sectional curvature constant equal to -1
and its fundamental group is two generated, infinitely
presented such that $\ker (\pi_1 (B(p,r)) \ra \pi_1(B(p,R)))$
is strictly smaller than $\ker (\pi_1 (B(p,r)) \ra
\pi_1(M,p))$ for all $r \geq 1, R\geq r$ \cite{BoMe},
\cite{Po}. So the universal covering of $M$ restricted to
$B(p,r)$ and the universal covering of $B(p,R)$ restricted to
$B(p,r)$ are always different for all $R\geq r$.

%
\end{ex}

We would like to
thank D. Cooper for bringing this and other similar examples to our
attention.

One consequence of this example is that we cannot
hope to study a ball in $\tilde{Y}$ just by
applying Theorem~\ref{limdelcov}.  This example demonstrates
that the lift of a ball $B(y,r)$ to $\tilde{Y}$ may not
be isometric to the limit of any relative delta cover no
matter how large we take $R$ and how small we take $\delta$.
The following theorem allows us to study regions in $\tilde{Y}$
using limits of relative delta covers with $R \to \infty$.

\begin{theo} \label{univisalim}
Given any $\tilde{x} \in \tilde{Y}$ and
$0<r< \infty$. Let $x=\pi(\tilde{x})\in Y$, and
$\bar{B}(x,r)$ the connected lift of $B(x,r)$ in $\tilde{Y}$
containing $\tilde{x}$.
Each space is given the intrinsic metric.
Then $(\bar{B}(x,r), \tilde{x})$ is the pointed Gromov
Hausdorff limit of stable relative delta
covers
$(\tilde{B}(x,r,R_i)^{\delta_i}=B(x,r,R_i)^{\delta_i},
\tilde{x}_i)$ where $\tilde{x}_i$ is a lift of $x_i$ and
$R_i>3r$ diverge to infinity. Furthermore these relative
delta covers all cover $\bar{B}(x,r)$. \end{theo}

Since each $B(x,r,R_i)^{\delta_i}$ is a GH limit of
$\tilde{B}(x_j, r_j, R_{i,j})^{\delta_i}$ by definition, we have the
following immediate corollary.

\begin{coro}\label{2univisalim}
In fact $(\bar{B}(x,r), \tilde{x})$ is the pointed Gromov
Hausdorff limit of a subsequence of relative $\delta$ covers
of $M_i^n$. Namely if we still index the subsequence as $k$,
then there are $r_k \to r$, $R_k \to \infty$, $B(x_k, R_k)
\subset M^n_k$ such that $(\tilde{B}(x_{k}, r_{k},
R_{k})^{\delta_k}, \tilde{x}_{k})$ converges to
$(\bar{B}(x,r), \tilde{x})$. \end{coro}
\Pf  We have $(\bar{B}(x,r), \tilde{x})$ is the pointed Gromov Hausdorff
limit of $(B(x,r,R_i)^{\delta_i}, \tilde{x}_i)$ for some $R_i
\geq 3r$ going to infinity and $\delta_i$ decreasing. So
$d_{GH} (B(x,r,R_i)^{\delta_i}, \bar{B}(x,r)) \leq
\epsilon_i$ for some $\epsilon_i \ra 0$. On the other hand
each $B(x,r,R_i)^{\delta_i}$ is a GH limit of $\tilde{B}(x_j,
r_j, R_{i,j})^{\delta_i}$. For $B(x,r,R_1)^{\delta_1}$,
choose $j_1$  sufficiently large such that $d_{GH}
(\tilde{B}(x_{j_1}, r_{j_1}, R_{i,j_1})^{\delta_1},
B(x,r,R_1)^{\delta_1}) \leq \epsilon_1$, for
$B(x,r,R_2)^{\delta_2}$, choose $j_2 \geq j_1$ sufficiently
large such that $d_{GH} (\tilde{B}(x_{j_2}, r_{j_2},
R_{i,j_2})^{\delta_2}, B(x,r,R_2)^{\delta_2}) \leq
\epsilon_2$, continue choosing a subsequence in this way, we
have $d_{GH} (\tilde{B}(x_{j_k}, r_{j_k},
R_{i,j_k})^{\delta_k},B(x,r,R_k)^{\delta_k}) \leq
\epsilon_k$. Rename the index $j_k$ as $k$ we have
$d_{GH}(\tilde{B}(x_{k}, r_{k}, R_{i,k})^{\delta_k},
\bar{B}(x,r)) \leq 2 \epsilon_k$. \qed

We now prove Theorem~\ref{univisalim} showing that regions in the
universal cover of the limit space are limits themselves.

\noindent {\bf Proof of Theorem~\ref{univisalim}:}
 First note that for $R>3r$,  $\bar{B}(x,r)$
can also be seen as the connected
lift of $B(x,r)$ in $\bar{B}(x,R)$ containing $\tilde{x}$.

By the definition of a cover, we know that for all $z\in B(x,R)$,
there is a $\delta_z>0$ such that $B_z(\delta_z)$ lifts isometrically
to $\bar{B}(x,R)$.  By compactness
of $B(x,R)$, there exists a $\delta_0>0$ such that any $B_z(\delta_0)$
lifts isometrically to the universal cover and thus to $\bar{B}(x,R)$.
By \cite[p81]{Sp} and Definition~\ref{defdel}, $\tilde{B}(x, R)^{\delta_0}$
must cover $\bar{B}(x,R)$.

So $\bar{B}(x,r)$ is covered by
$\tilde{B}(x,r,R)^{\delta_0}$.
Applying Theorem~\ref{samedelta}, we know that there exists $\delta_R>0$
such that the relative delta covers stabilize, so for all $\delta<\delta_R$,
we have
$\tilde{B}(x,r,R)^\delta
={B}(x,r,R)^{\delta}=\tilde{B}(x,r,R)^{\delta_{R}}$ all cover
$\bar{B}(x,r)$.

Now let us  take a sequence of $R_i$ diverging to infinity and let $\delta_i$
be chosen such that $\delta_i\le\delta_{R_i}$ and $\delta_i$ decrease.
Then by the stabilization
$\tilde{B}(x,r,R_i)^{\delta_{i+1}}=\tilde{B}(x,r,R_i)^{\delta_i}.$
Furthermore $\tilde{B}(x,r,R_i)^{\delta_{i+1}}$
covers $\tilde{B}(x,r,R_{i+1})^{\delta_{i+1}}$
because any curve that lifts to a closed curve in
$\tilde{B}(x,r,R_i)^{\delta_{i+1}}$ is homotopic in $B(x,R_i)$ to
a combination of curves of the form $\alpha \beta \alpha^{-1}$ with $\beta$
in a $\delta_{i+1}$ ball, and so this is also true if we allow the
homotopy in the larger ball $B(x, R_{i+1})$.

Thus we have
\be
\tilde{B}(x,r,R_i)^{\delta_i}
={B}(x,r,R_i)^{\delta_i}
\ra
\tilde{B}(x,r,R_{i+1})^{\delta_{i+1}}
={B}(x,r,R_{i+1})^{\delta_{i+1}}
\ra \cdots
\ra \bar{B}(x,r).
\ee
Let $\tilde{x}_i$ be a lift of $x_i$ to
$\tilde{B}(x,r,R_i)^{\delta_i}$.

{\em We claim that a subsequence of
$(\tilde{B}(x,r,R_i)^{\delta_i}, \tilde{x}_i)$
converges in the pointed Gromov Hausdorff sense to a limit space
$\bar{B}(x,r,\infty)$.}

To prove this we apply Gromov's Compactness Theorem
\cite{Pe2}[Page 280, Lemma 1.9][Gr].  Given any $\rho>0$,
$\epsilon>0$ we must show that for $i$ sufficiently large
there is a uniform bound on the number of disjoint balls of
radius $\epsilon$ contained a ball $B(\tilde{x}_i, \rho)$.
We can just take $i$ large enough that $R_i>4\rho$, and then
apply Proposition~\ref{limmeas} using the fact that
$\tilde{B}(x,r,R_i)^{\delta_i}=B(x,r,R_i)^{\delta_i}$
with renormalized limit measures.

{\em We claim $\bar{B}(x,r,\infty)$ is a cover of $\bar{B}(x,r)$.}

By Lemma~\ref{localisom} we know that
the covering map $\pi_1:\tilde{B}(x,r,R_1)^{\delta_1}\to \bar{B}(x,r)$
is an isometry on balls of radius $\delta_1/3$.  Since,
we have covers
$f_i:\tilde{B}(x,r,R_1)^{\delta_1} \to \tilde{B}(x,r,R_i)^{\delta_i}$
and $\pi_i:\tilde{B}(x,r,R_i)^{\delta_i}\to \bar{B}(x,r)$, each
$\pi_i$ must preserve the same isometry of balls of radius $\delta_1/3$.
Furthermore the $\pi_i$ are uniformly equicontinuous, so by a generalized
Arzela-Ascoli theorem (see e.g. \cite[Page 279, Lemma 1.8]{Pe2})
a subsequence converges to a continuous function
$\pi_\infty:\bar{B}(x,r,\infty) \to \bar{B}(x,r)$
which is an isometry on balls of radius less than $\delta_1/3$.
Thus $\pi_\infty$ is a covering map.

{\em We claim $\bar{B}(x,r,\infty)$ is isometric to $\bar{B}(x,r)$.}

Since $\bar{B}(x,r, \infty)$ is a cover of $\bar{B}(x,r)$
then we need only show that
all loops in $B(x,r)$ which lift
non closed to $\bar{B}(x, r, \infty)$ also lift nonclosed to
$\bar{B}(x,r)$ \cite [Page 78, Lemma 9] {Sp}.
But if a loop $C$ lifts nonclosed to $\bar{B}(x,r,\infty)$
then by the Hausdorff approximation
it lifts nonclosed  to all $B(x,r, R_i)^{\delta_i}$ for $i$ sufficiently large
depending on $C$.
So $C$ is not homotopic in $B(x,R_i)$ to a combination of loops
$\alpha \beta \alpha^{-1}$
with $\beta$ in a ball of radius $\delta_i$.  In particular $C$ is not
contractible in $B(x, R_i)$.  But this is true for $R_i$ arbitrarily large.
Thus $C$ must not be
contractible for if it were contractible there would be a homotopy
$H:I \times I \to Y$
but $im(H)$ would be compact and fit in some $B(x, R_i).$

\qed

Our first application of Theorem~\ref{univisalim}, will be to
define a renormalized limit measure on $\tilde{Y}$ and
prove that it satisfies the Bishop Gromov Volume comparison
globally.  We begin with a lemma.

 \begin{lem} \label{natmeas}
Given a
length space $X_1$ with a measure $\mu_1$ and  a
covering $\pi_{12}:X_2 \to X_1$, there is a natural lifted measure
$\mu_2=\mu_{\pi_{12}}$
 on $X_2$ so that the covering
map $\pi_{12}$ is locally measure preserving and globally measure
nonincreasing.

Furthermore, if $X_3\to X_2 \to X_1$ are all coverings
then $\mu_3=\mu_{\pi_{13}}$ where $\pi_{13}:X_3\to X_1$,
agrees with $\mu_{\pi_{23}}$, the natural lifted measure
from $X_2$ to $X_3$ of $\mu_{\pi_{12}}$.
 \end{lem}

\Pf We can define a function $\Phi$ on small open balls in
$X_2$ which are projected isometrically to $X_1$
by pulling back the measure on $X_1$.
Now for any closed subset $A \subset X_2$, we define
$\bar{\mu} (A)$ by the Caratheodory's construction (see
\ref{rcarath1} and \ref{rcarath2}).   Since the measure
on $Y$ was also constructed using Caratheodory based on the
same $\Phi$ \cite{ChCo2}, the covering map $\pi$ is locally
measure preserving.  This gives a Borel measure on
$X_2$.  Since the $\Phi$ will agree locally on $X_3$ as well,
it defines the same lift in two steps or in one.

  One can use packing arguments to see that
$\pi$ is measure nonincreasing.  
\qed

 \begin{theo} \label{VOL}
 The
universal cover $\tilde{Y}$ with the natural lifted measure
$\tilde{\mu}=\mu_{\pi}$,
of the renormalized limit measure $\mu_Y$ on $Y$,
satisfies Bishop-Gromov volume comparision
globally,
\be
\frac
{\tilde{\mu}(B(\tilde{x},r_1))}{\tilde{\mu}(B(\tilde{x},r_2))
} \ge V(n,H,r_1)/V(n,H,r_2) \ \ \forall r_1\le r_2,\
\tilde{x} \in \tilde{Y}.  \ee
\end{theo}

\Pf
By Theorem~\ref{univisalim}, we know that for any fixed $r$,
$(\bar{B}(y,r), \tilde{y})$ is the pointed Gromov Hausdorff
limit of stable relative delta
covers $(\tilde{B}(y,r,R_i)^{\delta_i}=B(y,r,R_i)^{\delta_i},
\tilde{y}_i)$ where $\tilde{y}_i$ is a lift of $y$ and $R_i$
diverge to infinity. Furthermore these relative delta covers
all cover $\bar{B}(y,r)$. In particular by
Lemma~\ref{localisom}, all these relative delta covers act as
isometries on balls of radius $\delta_1/4$.  Note that
$\delta_1$ depends on $r$ so we will let
$\delta_r=\delta_1/2$.

Now each $(B(y,r,R_i)^{\delta_i}, \tilde{y}_i)$
is a limit of $(\tilde{B}(p_j,r_j,R_{i,j})^{\delta_i},
\tilde{p}_j)$, so for each $i$ we can take $j_i$ sufficiently
large that \be
d_{GH}((B(y,r,R_i)^{\delta_i}, \tilde{y}_i),
(\tilde{B}(p_{j_i},r_{j_i},R_{i,j_i})^{\delta_i},
\tilde{p}_{j_i}))<\epsilon_i \ee
where $\epsilon_i <\min\{ \delta_r/100, 1/i\}$.
Thus $(\tilde{B}(p_{j_i},r_{j_i},R_{i,j_i})^{\delta_i},
\tilde{p}_{j_i})$ converge to $\bar{B}(y,r)$ as well.

Furthermore, we claim that
$\pi_i:\tilde{B}(p_{j_i},r_{j_i},R_{i,j_i})^{\delta_i}\to B(p_i,r_{j_i})$
is an isometry on balls of radius $\delta_r/4$.  If this were not true
then there would be a pair of lifts $\tilde{p}_i$ and $\tilde{p}_i'$
of $p_i$ such that $d(\tilde{p}_i, \tilde{p}_i')< \delta_r/2$.
By the convergence of the covering maps in Theorem~\ref{limdelcov}, 
and the fact that we are within $\epsilon_i$ of the limit, this implies that
there are two lifts $\tilde{y}_i$ and $\tilde{y}_i'$ of $y$ in 
$B(y,r_i,R_i)^{\delta_i}$ such that
$d(\tilde{y}_i, \tilde{y}_i')< \delta_r/2+2\epsilon_i<\delta_r$.
This contradicts the isometry of $\delta_r/2$ balls on these relative
delta covers.

Thus we can apply Proposition~\ref{limmeas} with $R=\infty$, and 
$\delta=\delta_r$, to get a renormalized limit measure
$\mu_r$ defined on $\bar{B}(y,r)\subset \tilde{Y}$
which satisfies the Bishop Gromov Volume Comparison for pairs of balls,
$B(x,r_1), B(x,r_2)$ such that the radius of the inner ball 
satisfies $r_1<r-d_{Y}(\pi(x), y)$.  Here the balls are measured
using $d_{\bar{B}(y,r)}$ and the outer radius can be arbitrarily
large.  To get Bishop Gromov for balls defined using $d_{\tilde{Y}}$ we apply
Lemma~\ref{distance}, and restrict the outer radius 
$r_2<r/3-d_{Y}(\pi(x),y)$.

Furthermore, by Corollary~\ref{agree} and
Lemma~\ref{natmeas}, this renormalized limit measure agrees
with the lifted measure $\mu_{\pi}$ as follows:
\be
\mu_r(S)= \lambda_r \,\mu_{\pi}(S),
\ee
where
\be
\lambda_r=\lim_{i\to\infty}
\frac{Vol(B(p_i,1) \subset M_i)}{Vol(B(\tilde{p}_i, r/10)\subset
      \tilde{B}(p_{j_i},r_{j_i},R_{i,j_i})^{\delta_i} )}
\in \left[\frac{ V(n,H,\delta/2)}{V(n,H,r/10)},
\frac{V(n,H,1)}{V(n,H,\delta/2)}\right].
\ee
 Since the Bishop Gromov Volume
Comparison is a ratio, it also holds with respect to the
lifted measure $\mu_{\pi}$ for balls measured with
$d_{\tilde{Y}}$ of radius bounded by $r$ as above.  However,
this is true for any $r>0$, so Bishop Gromov holds for balls
of all sizes and locations. \qed

We also get a splitting theorem on $Y$ if the sequence of manifolds
have Ricci curvature converging towards a nonnegative lower bound.

\begin{theo} \label{SPLIT}
Let $(M_i^n,p_i)$ be a sequence of manifolds with $\Ric_{M_i}
\geq -(n-1)\epsilon_i, \epsilon_i \ra 0$ and converges to
$(Y,y)$ in the pointed Gromov-Hausdorff sense. By
Theorem~\ref{mainthm},  the universal cover of $Y$, $\tilde{Y}$,
exists. If $\tilde{Y}$ contains a line, then $\tilde{Y}$
splits isometrically, $\tilde{Y} = {\mathbb R} \times
\bar{X}$.
\end{theo}

\Pf Let us assume the line of $\tilde{Y}$ passes
through $\tilde{x} \in \tilde{Y}$. Given any $R> 0, L >
1$, let $B(\tilde{x}, 2LR)$ be a ball in $\tilde{Y}$. Then
there is a line segment $\tilde{\gamma}: (-2LR, 2LR) \ra
B(\tilde{x}, 2LR)$ with $\tilde{\gamma}(0) = \tilde{x}$.
Let $x = \pi(\tilde{x})$ and $\bar{B}(x,6LR)$ be the
connected lift of $B(x,6LR)\subset Y$ in $\tilde{Y}$
containing $\tilde{x}$. Then $B(\tilde{x},2LR) \subset
\tilde{Y}$ are the same as $B(\tilde{x},2LR) \subset
\bar{B}(x,6LR)$. By Corollary~\ref{2univisalim}
$(\bar{B}(x,6LR), \tilde{x})$ is the pointed Gromov-Hausdorff
limit of a sequence of relative $\delta$-covers of manifolds
with $\Ric \geq -(n-1)\epsilon_{k}$, namely
$(\tilde{B}(x_{k}, r_{k}, R_{k})^{\delta_k},
\tilde{x}_{j_i})$, where $r_{k} \to 6LR$, $R_{k}\to
\infty$. Let $q_{k}^+, q_{k}^- \in \tilde{B}(x_{k},
r_{k}, R_{k})^{\delta_k}$ be the $\epsilon_{k}$
Hausdorff images of the points
$\tilde{\gamma}(2LR),\tilde{\gamma}(-2LR)$. Now given any
$\epsilon > 0$, choose $L = L(n,\epsilon)$ the constant in
\cite[Proposition 6.2]{ChCo1}, then $\min
(d(x_{k},q_{k}^+), d(x_{k},q_{k}^-)) \geq LR$, the
excess at $x_{k}$ with respect to $q_{k}^+, q_{k}^-$,
$E(x_{k}) \leq 3\epsilon_{k}$. Since $\epsilon_{k} \ra
0$ as $k \ra \infty$. We have for all $k$ big,
$E(x_{k}) \leq \tau R$ and $\Ric \leq -(n-1) \tau R^{-2}$.
By \cite[Theorem 6.62]{ChCo1}, the ball $B(x_{k},R) \subset
\tilde{B}(x_{k}, r_{k}, R_{k})^{\delta_k}$ is $\epsilon
R$ Gromov-Hausdorff close to an $R$-ball in ${\mathbb R}
\times X_{k}$, for some metric space $X_{k}$. Therefore
the limit ball $B(\tilde{x},R)$ is isometric to an $R$-ball
in ${\mathbb R} \times X$, for some metric space $X$. This is
true for any $R$ ball in $\tilde{Y}$. Hence $\tilde{Y}$
splits globally. \qed

Using Theorem~\ref{VOL} we can easily extend several
results about manifolds with nonnegative Ricci curvature to
limit spaces.

First we can extend Milnor's result \cite{Mi} about
fundamental groups of polynomial growth to the revised
fundamental groups of limit spaces.

 \begin{coro} \label{cormil}
 If $(Y,y)$ is a pointed Gromov-Hausdorff limit of a sequence
of complete manifolds $(M_i^n,p_i)$ with $\Ric_{M_i} \ge
-(n-1)\epsilon_i, \epsilon_i \ra 0$,
then any finitely generated subgroup of the revised
fundamental group of $Y$, $\bar{\pi}_1(Y)$, is of polynomial
growth of degree at most $n$.
\end{coro}

We can also extend Anderson's Theorems from \cite{An2},
regarding volume growth and the revised fundamental group.
In particular

 \begin{coro} \label{corand}
 If $(Y,y)$ is a pointed Gromov-Hausdorff limit of a sequence
of complete manifolds $(M_i^n,p_i)$ with $\Ric_{M_i} \ge 0$,
and if it has Euclidean measure growth 
$\liminf_{r\to\infty} \mu(B_{\tilde{y}}(r))/r^n=C>0$
then the revised fundamental group is finite
and $|\bar{\pi}_1(Y,y)|\le \omega_n/C$.
\end{coro}

We will say that a length space $Y$ has the loops
 to infinity property if given any element $g$ of
 the revised fundamental group of $Y$ based at $y$ and
 given any compact set $K$ in $Y$, $g$ has a representative
 element of the form $\gamma \circ C \circ \gamma^{-1}$
 where $C$ is a loop in $M\setminus K$ and $\gamma$ is a minimal curve running
 from $y$ to $M\setminus K$.

Using this definition it is easy to imitate the beginning of the proof in [So2]
to obtain the following:

 \begin{coro} \label{corsor}
 If $M_i$ have $Ricci(M_i) \ge -\epsilon_i$
 with $\epsilon_i$ decreasing to 0, and if
 $Y=GH\lim_{i\to\infty} M_i$ then
 either $Y$ has the loops to infinity
 property or the universal cover of $Y$
 splits isometrically.
 \end{coro}

 Finally we close with a theorem relating the local fundamental groups
 of the $M_i$ to the revised fundamental groups in Y.

\begin{prop} \label{localsurj}
$G(p_i,r_i,R_i)$ maps surjectively to
$im(\bar{\pi}_1(B(y,r),y) \ra \bar{\pi}_1(Y,y))$ for $i \ge
N_{r,R,\delta}$, where $G(p_i,r_i,R_i)$ is the deck
transformations of the connected lift of $B(p_i,r_i)$ in the
universal cover of $B(p_i,R_i)$.
 \end{prop}

Recall Example~\ref{cylinder} in which we showed that we can not hope
for  surjectivity when we do not restrict our attention to
balls.  Recall  also Example~\ref{cooper} which
shows that even if $M_i=Y$, we  can have a large kernal.
Also $M_i$ to $Y$ figure eights to a cylinder  demonstrate
that the kernal may be free even without collapsing.  This is
in strong contrast with the compact case in which the kernal
is finite  in the noncollapsed case.

 \Pf  By
Theorem~\ref{samedelta} for any $r < R$,
$\tilde{B}(y,r,R)^\delta$ stabilizes for
$\delta<\delta_{y,r,R}$.  So $G(y,r,R,\delta)$ also
stabilizes for $\delta<\delta_{y,r,R}$.  Call the stable
group $G(y,r,R)$.

By Corollary~\ref{GG} we have surjective maps from
$G(p_i,r_i,R_i, \delta_1)$ to this stabilized group
$G(y,r,R)$ for all $\delta_1<\delta_{y,r,R}$ and $i$ large.
Now there are natural surjective maps from
\be
G(p_i,r_i,R_i)=im(\pi_1(B(p_i,r_i),p_i)\ra
\pi_1(B(p_i,R_i),p_i))
\ee 
onto
$G(p_i,r_i,R_i,\delta_1)$.
So we have surjective maps from  $G(p_i,r_i,R_i)$ to
$G(y,r,R)$ for $i \ge N_{r,R,\delta}$.

However $G(y,r,R)$ maps surjectively onto
$im(\bar{\pi}_1(B(y,r),y)\ra \bar{\pi}_1(Y,y))$
because $\tilde{B}(y,r,R)^{\delta_{y,r,R}}$ covers the
connected lift of $B(y,r)$ in $\tilde{Y}$.

Thus $G(p_i,r_i,R_i)$ maps onto $im(\bar{\pi}_1(B(y,r),y)\ra
\bar{\pi}_1(Y,y))$  for $i \ge N_{r,R,\delta}$.
\qed

  \sect{Appendix}

Once and for all we show the two possible definitions of Gromov Hausdorff
convergence on noncompact spaces are identical.

\begin{defn}\label{GrRest}
We say $(M_i,p_i)$ converges in the pointed Gromov Hausdorff
sense to $(Y,y)$ if for all $R>0$, $(B(p_i,R), d_{M_i})$
converges to $(B(y,R), d_Y)$ in the Gromov Hausdorff sense.
\end{defn}

This definition works well because if $M_i$ happen
to have a uniform upper bound on diameter and $M_i$ converge
to $Y$, then there exist $p_i$ and $y$ such that
$(M_i,p_i)$ converges in the pointed Gromov Hausdorff
sense to $(Y,y)$. See \cite[Page 279]{Pe2}.

On the other hand, Gromov emphasizes the importance of length spaces
in his text.  For our purposes it is essential to use length spaces.
Thus there is another possible definition.

\begin{defn}\label{GrInt}
We say $(M_i,p_i)$ converges in the pointed Gromov Hausdorff
sense to $(Y,y)$ if for all $R>0$, there exists $R_i \to R$
such that $(B(p_i,R_i), d_{B(p_i,R_i)})$
converges to $(B(y,R), d_{B(y,R)})$ in the Gromov Hausdorff sense.
\end{defn}

Note that if $M_i$ converges to $X$ in Gromov-Hausdorff
topology and $p_i \in M_i$ converges to $x \in X$, it may not true
that $B(p_i,R)$ in $M_i$ with intrinsic metric converges to
$B(x,R)$ in $X$ with intrinsic metric. For example, let $M_i$ be
circles of radius $1+(1/i)$ converges to $X$, the circle with
radius $1$. Then the balls of radius $\pi$ in $M_i$ with intrinsic
metric is the intervals $[0,2\pi]$ while the ball of radius
$\pi$ in $X$ is the unit circle. Therefore in
above one needs $R_i$'s not just $R$.

\begin{lem} \label{grdefsame}
Definition~\ref{GrRest} and Definition~\ref{GrInt} are
equivalent.
\end{lem}

Before proving this statement, we make the following definition.

\begin{defn} \label{HausApprx}
Given two metric spaces $X,Y$, a map $\varphi: X \ra Y$ is said to be an
$\epsilon$-Hausdorff approximation if the following conditions are satisfied.

1) The $\epsilon$-neighborhood of $\varphi (X)$ in $Y$ is equal to $Y$.

2) For each $x_1, x_2$ in $X$, we have
\[
| d_X(x_1,x_2) - d_Y(\varphi (x_1), \varphi (x_2))| \leq \epsilon.
\]
\end{defn}

We say metrics spaces $X_i$ converges to $X$ in Gromov-Hausdorff
topology if there are
$\epsilon_i$-Hausdorff approximations from $X_i$ to $X$ and
$\epsilon_i \ra 0$ as $i \ra \infty$.

\Pf
We first show that if $(M_i,p_i)$ converges in the pointed Gromov Hausdorff
sense to $(Y,y)$ according to Definition~\ref{GrInt} then it
converges according to Definition~\ref{GrRest} as well.

For all $R>0$, we know  there exists $\rho_i \to 3R$
such that $(B(p_i,\rho_i), d_{B(p_i,\rho_i)})$
converges to $(B(y,3R), d_{B(y,3R)})$ in the Gromov Hausdorff sense.
Thus $(B(p_i,R), d_{B(p_i,\rho_i)})$
converges to $(B(y,R), d_{B(y,3R)})$ in the Gromov Hausdorff sense
as well.  However by Lemma~\ref{distance},  $(B(p_i,R),
d_{B(p_i,\rho_i)})$ is isometric to  $(B(p_i,R), d_{M_i})$ and
$(B(y,R), d_{B(y,3R)})$ is isometric to $(B(y,R), d_{Y})$  and we
are done.

The other direction is less trivial.

Given $(M_i,p_i)$ converges in the pointed Gromov Hausdorff
sense to $(Y,y)$ according to Definition~\ref{GrRest}.  Fix $R>0$.
We must construct $R_i\to R$.

First there exists $\epsilon_i \to 0$ such that
\be
d_{GH}\left((B(p_i,3R), d_{M_i}),(B(y,3R), d_Y) \right) <\epsilon_i
\ee
so there exists
\be
f_i:B(y,3R) \ra B(p_i,3R) \textrm{ where } f_i(y)=p_i
\ee
which is $\epsilon_i$-Hausdorff approximation with restricted metrics.

Now $f_i:B(y,R) \to B(p_i, R+\epsilon_i)$ and we would like to
show that it is almost distance preserving and almost onto
with the intrinsic distances.  However, we will not be able to
do so without adding a little extra space.  So we will look at
\be
f_i:B(y,R) \to B(p_i, R+\epsilon_i+\delta_i)
\ee
where $\delta_i=\sqrt{\epsilon_i}$.

We will show that $\forall \delta>0$ there exists $i$ sufficiently
large such that $f_i: B(y,R) \to B(p_i, R+\epsilon_i+\delta_i)$ are
$\delta$-Hausdorff approximation with the intrinsic distances.
Note that $f_i$ could be far from being such a map if $i$ is not
taken sufficiently large.  Take $Y=M_i =$ cylinders
such that $B(y,R)$ does not have cut points but
$B(p_i, R+\epsilon_i+\delta_i)$ does.  Note that for $i$ sufficiently
large $B(p_i, R+\epsilon_i+\delta_i)$ doesn't have cut points anymore.

We begin by showing it is almost distance preserving.
 For all $a,b\in B(y,R)$ there is a curve $C\in B(y,R)$ from
$a$ to $b$, such that $L(C)=d_{B(y,R)}(a,b) \le 2R$.  We parametrize
$C$ by arclength and take $t_0=0, t_j=t_{j-1}+\delta_i$ and
$t_N\leq t_{N-1}+\delta_i$ so that $C(t_0)=a, C(t_N)=b$,  
\be
d_Y(C(t_j),C(t_{j+1})) \le d_{B(y,R)}(C(t_j),C(t_{j+1}))\le \delta_i
\ee
and
\be
L(C)-\delta_i \le (N-1)\delta_i \le L(C).
\ee

Now $f_i(a)$ and $f_i(b)$ are joined by a curve $C_i$
which is created by joining the points $f_i(t_j)$ by
minimal curves in $M_i$.  Note that then 
$C_i \subset B(p_i, R+\epsilon_i+\frac{\epsilon_i+\delta_i}{2})$,
and
\begin{eqnarray}
L(C_i) &=&\sum_{j=1}^N d_{M_i}\left(f_i(C(t_{j-1})), f_i(C(t_j))\right) \\
&\le& \sum_{j=1}^N \left( d_{Y}\left(C(t_{j-1}), C(t_j)\right) +\epsilon_i 
\right) \\
&\le& N(\delta_i+\epsilon_i) = N\delta_i(1+\delta_i) \\
&\le & (L(C)+ \delta_i)(1+\delta_i) \le L(C) +(2R+1)\delta_i+\epsilon_i.
\end{eqnarray}

Thus 
\be \label{eqnadp1}
d_{B(p_i,R+\epsilon_i+\delta_i)}(f_i(a), f_i(b)) \le
d_{B(y,R)}(a,b)+(2R+1)\delta_i+\epsilon_i.
\ee

On the other hand,
$d_{B(p_i,R+\epsilon_i+\delta_i)}(f_i(a),f_i(b))=L(\sigma_i)= L_i$
for some curve $\sigma_i:[0,L_i]\to B(p_i,R+\epsilon_i+\delta_i)$ connecting 
$f_i(a),f_i(b)$.
Let $L_\infty=\liminf_{i\to\infty} L_i$.  Then there is a subsequence
converging to $L_\infty$.  Since $(B(p_i,R+\epsilon_i+\delta_i), d_{M_i})$ 
converges to $(B(y,R), d_Y)$, by the generalized Arzela-Ascoli
theorem \cite{GP}, \cite[page 279, Lemma 1.8]{Pe2}, there is 
a subsequence such that $\sigma_i$ converges to
$\sigma_\infty:[0,L_\infty] \to B(y, R)$ and
$\sigma_\infty(0)=a$ and $\sigma_\infty(L_\infty)=b$.

Thus 
$d_{B(y,R)}(a,b)\le L(\sigma_\infty)=\lim_{i\to\infty}L(\sigma_i)=L_\infty$.
Since $L_\infty=\liminf_{i\to\infty} L_i$, 
for all $\delta >0$ there exists $N_{a,b,\delta}$ sufficiently large such that
\be \label{eqnab}
d_{B(y,R)}(a,b)\le d_{B(p_i,R+\epsilon_i+\delta_i)}(f_i(a),f_i(b))+\delta \ \ \ 
\ \forall \ i\ge N_{a,b,\delta}.
\ee
However we need a uniform estimate for $N$ not depending on $a$ and $b$
to say that $f_i$ is almost distance preserving.

We assume on the contrary that there exists $\delta>0$
and a subsequence $i\to \infty$ and points $a_i,b_i\in B(y,R)$ such that
\be \label{eqnadp2}
d_{B(y,R)}(a_i,b_i)\ge 
d_{B(p_{i},R+\epsilon_{i}+\delta_{i})}(f_i(a_i),f_i(b_i))+\delta.
\ee

Then a subsequence of $a_i$ and of $b_i$ converge to points $a$ and $b$
in $B(y,R)$ for which (\ref{eqnab}) holds.  So take $i\ge N_{a,b,\delta/10}$
and also large enough that $(2R+1)\delta_i+\epsilon_i<\delta/2$,
$d_{B(y,R)}(a,a_i)<\delta/10$ and $d_{B(y,R)}(b,b_i)<\delta/10$.  
Then we can apply (\ref{eqnab}) and (\ref{eqnadp1}) to get
\begin{eqnarray}
d_{B(y,R)}(a_i,b_i) &\le&  2(\delta/10) + d_{B(y,R)}(a,b) \\
&\le& 3(\delta/10)+d_{B(p_i,R+\epsilon_i+\delta_i)}(f_i(a),f_i(b)) \\ 
&\le & 3(\delta/10)+ d_{B(y,R)}(a,b)+ (2R+1)\delta_i+\epsilon_i 
<d_{B(y,R)}(a,b)+\delta,
\end{eqnarray}
which contradicts (\ref{eqnadp2}).

Thus (\ref{eqnab}) holds without the dependence on $a$ and $b$.  
Combining this with (\ref{eqnadp1}) we have,
$\forall \delta>0, \exists N_\delta$ such that $f_i$ is $\delta$ 
almost distance preserving for all $i\ge N_\delta$.

Now to prove that $f_i$ is $\delta$ almost onto, we note that
\be
\forall q\in B(p_i,R+\epsilon_i+\delta_i)\subset B(p_i, 3R)  \ \ 
\exists z'_q\in B(y,3R) \ s.t. \
d_{M_i}(f_i(z'_q), q) <\epsilon_i.
\ee

If $q\in B(p_i, R-2\epsilon_i)$ then
\begin{eqnarray}
d_Y(z'_q,y)&\le&d_{M_i}(f_i(z'),p_i)+\epsilon_i \\
&\le& d_{M_i}(q, p_i) + d_{M_i}(f_i(z'), q)+\epsilon_i
<R-2\epsilon_i+2\epsilon_i
\end{eqnarray}
So  $z'_q \in B(y,R)$.
Furthermore the minimal curve from $f_i(z'_q)$ to $q$ must
be in $B(p_i,R-\epsilon_i)$
so
\be
d_{ B(p_i,R+\epsilon_i+\delta_i)}(z'_q,q)=d_{M_i}(z'_q,q)<\epsilon_i.
\ee

If $q\in B(p_i,R+\epsilon_i+\delta_i)$, then let $\bar{q}$
be the first point on a minimal geodesic joining $q$ to $p_i$
which is in $B(p_i, R-2\epsilon_i)$.
So $d_{ B(p_i,R+\epsilon_i+\delta_i)}(q,\bar{q})<3\epsilon_i+\delta_i$.
Then let $z_q=z'_{\bar{q}}$, so
\begin{eqnarray}
d_{ B(p_i,R+\epsilon_i+\delta_i)}(z_q,q)
&\le& d_{ B(p_i,R+\epsilon_i+\delta_i)}(q, \bar{q})
+d_{ B(p_i,R+\epsilon_i+\delta_i)}(z_q,\bar{q})  \\
&<&3\epsilon_i+\delta_i+\epsilon_i.
\end{eqnarray}

Thus $f_i$ is $4\epsilon_i+\delta_i$ almost onto with respect to
the intrinsic distances.

Since $4\epsilon_i+\delta_i$ converges to $0$, we are done.
\qed

\begin{coro}  \label{R-r}
If $(M_i,p_i)$ converges in the pointed Gromov Hausdorff
sense to $(Y,y)$, then for any $0<r<R$, there exist
$r_i \ra r, R_i \ra R$, such
that for all $\delta>0$ there exists $N_\delta(r,R)$ and maps $f_i$ such
that for all $i \geq N_\delta(r,R)$ the maps $f_i$ are $\delta$-Hausdorff
approximations $f_i : B(y,R) \to B(p_i, R_i)$ with
respect to intrinsic metrics,
and their restrictions $f_i:B(y,r)\to B(p_i,r_i)$ are also $\delta$-Hausdorff
approximation with respect to the intrinsic distances
on these smaller balls.
\end{coro}

Department of Mathematics and Computer Science,

Lehman College, City University of New York,

Bronx, NY 10468

sormani@g230.lehman.cuny.edu

Department of Mathematics, 

University of California, 

Santa Barbara, CA 93106 

wei@math.ucsb.edu

\end{document}